\documentclass[11pt,reqno]{amsart}
\usepackage{amssymb,latexsym}
\usepackage{graphicx}
\usepackage{comment}
\usepackage{url}		
\usepackage{color}
\usepackage{hyperref}
\usepackage{url}
\usepackage{breakurl}
\newcommand{\bburl}[1]{\textcolor{blue}{\url{#1}}}

\calclayout

\makeatletter
\newcommand{\monthyear}[1]{%
	\def\@monthyear{\uppercase{#1}}}
\newcommand{\volnumber}[1]{%
\def\@volnumber{\uppercase{#1}}}
\AtBeginDocument{%
\def\ps@plain{\ps@empty
\def\@oddfoot{\@monthyear \hfil \thepage}%
\def\@evenfoot{\thepage \hfil \@volnumber}}
\def\ps@firstpage{\ps@plain}
\def\ps@headings{\ps@empty
\def\@evenhead{%
	\setTrue{runhead}%
	\def\thanks{\protect\thanks@warning}%
	\uppercase{}\hfil}%
\def\@oddhead{%
	\setTrue{runhead}%
	\def\thanks{\protect\thanks@warning}%
	\hfill\uppercase{}}%
\let\@mkboth\markboth
\def\@evenfoot{%
	\thepage \hfil \@volnumber}%
\def\@oddfoot{%
	\@monthyear \hfil \thepage}%
	}%
	\footskip=25pt
	\pagestyle{headings}%
}
\makeatother

\usepackage{amsthm}
\usepackage{booktabs}

\theoremstyle{plain}
\numberwithin{equation}{section}
\newtheorem{thm}{Theorem}[section]
\newtheorem{theorem}[thm]{Theorem}
\newtheorem{lemma}[thm]{Lemma}

\newtheorem{remark}[thm]{Remark}

\newtheorem{conjecture}{Conjecture}

\begin{document}
\monthyear{}
\volnumber{}
\setcounter{page}{1}

\title{Egg Drop Problems: They are all They are cracked up to be!}

\author{Xiangwen Cao, Zongyun Chen, Steven J. Miller}

\address{Building B, Kangyi Community, Yunlong District, Xuzhou, Jiangsu Province, China}
\email{caoxiangwen2008@outlook.com}

\address{193/198, Soi 9, Koolpunt Ville 7, Moo 4, Muang Chiangmai District, Chiangmai Province, Thailand}
\email{marinachenzongyun@gmail.com}

\address{Department of Mathematics\\
Williams College\\
Williamstown, MA\\}
\email{sjm1@williams.edu}

\thanks{This paper is an outgrowth of a talk given by the third named author and attended by the first and second named authors at the Math League International Summer Tournament in Summer 2025 in Trenton, NJ. It is a pleasure to thank the organizers and participants for creating such a lively atmosphere, especially John Cui, Dan Flegler, John Hagen, Rui Hu, and Adam Raichel.}

\maketitle

\begin{abstract}  We illustrate how to invite and excite students about research by exploring higher-dimensional generalizations of the classical egg drop problem, in which the goal is to locate a critical breaking point using the fewest number of trials. Beginning with the one-dimensional case, we prove that with $k$ eggs and $N$ floors, the minimal number of drops in the worst case satisfies $P_1(k) \leq \lceil k N^{1/k} \rceil$. We then extend the recursive algorithm to two and three dimensions, proving similar formulas: $P_2(k) \leq \lceil (k-1)(M+N)^{1/(k-1)} \rceil $ in 2D and $P_3(k) \leq \lceil (k-2)(L+M+N)^{1/(k-2)} \rceil$ in 3D, and conjecture a general formula for the $d$-dimensional case. Beyond the critical point problems, we then study the critical line problems, where the breaking condition occurs along $x+y=V$ (with slope $-1$) or, more generally, $\alpha x+\beta y=V$ (with the slope of the line unknown). We discuss how one frequently has to pivot from the original problem, which is intractable, to something that can be solved; in our case, using induction and recursion, two standard proof techniques.
\end{abstract}

\section{Motivation}

There are two complementary goals of this work. The first is to describe a successful approach to help students transition from taking straightforward mathematics classes to doing research. In particular: How does one find good questions to ask, and when those prove too difficult how can one pivot to find something that is still of interest and solvable? The second, of course, is to ask and solve such questions!

\begin{figure}[h]
\begin{center}
	\scalebox{.15}{\includegraphics{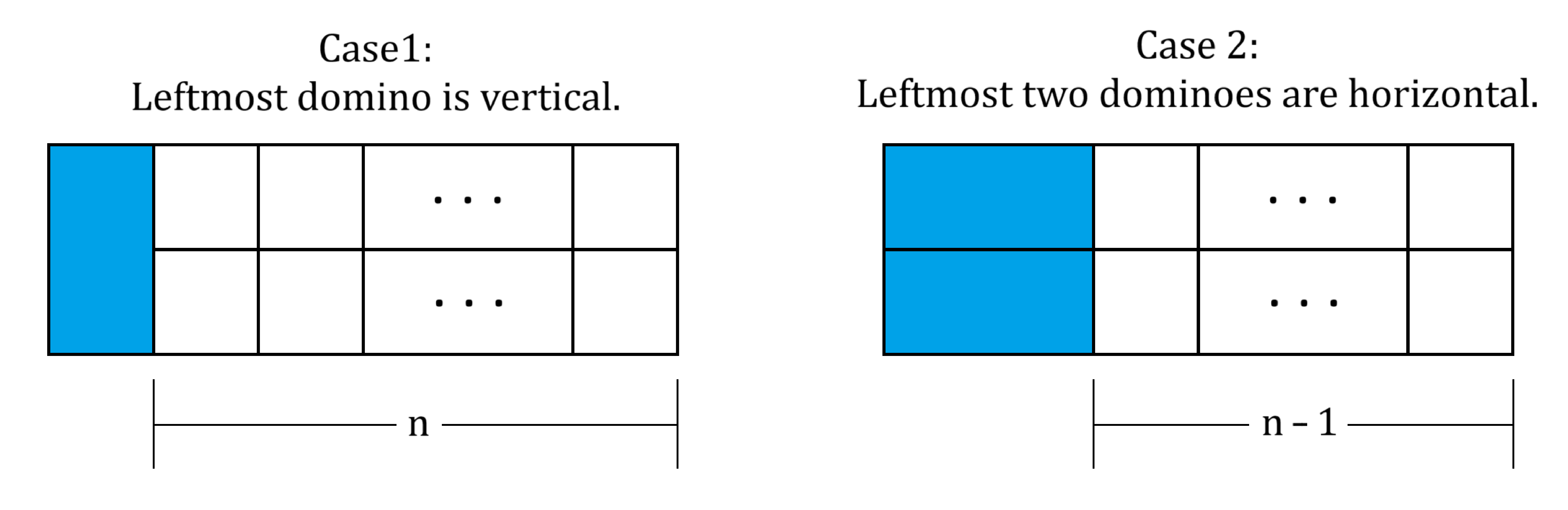}}
	\caption{\label{fig:Dominoes} Solution to the recursive tiling problem.}
	\end{center}\end{figure}

We investigate the famous egg drop problem and explore several generalizations. These problems share a common theme: we first solve a problem with a smaller value of the key parameter, and derive a recurrence relation to express the solution when we increase by 1 to the earlier case. This is an incredibly powerful technique, and many common problems are solved this way. For example, imagine we have a $n \times 2$ box and want to cover it with $2 \times 1$ dominoes that can be oriented either vertically or horizontally. If we let $T_n$ be the number of such tilings we quickly see that $T_1 = 1$, $T_2 = 2$ and $$T_{n+1} \ = \ T_{n} + T_{n-1}$$ (if the far left is covered by a vertically oriented domino we are left with a $n \times 2$ box to cover, while otherwise if we start with two horizontally oriented dominoes we then have a $(n-1) \times 2$ box to cover; clearly the two cases are distinct as one starts with a vertical and the other with two horizontal tiles, as illustrated in Figure \ref{fig:Dominoes}).
	Thus if we know how many tilings there are for boxes of length at most $n$, the recurrence immediately gives us the number for $n+1$ (recognizing this as the Fibonacci recurrence, we see the number of tilings are the Fibonacci numbers).
	
	Our proofs also provide an opportunity to highlight a good aspect of mathematical writing: if you have similar results used throughout a paper, it's best to prove a general result with free parameters, and then specialize as needed for the different situations. For a good example of this, see Lemma 2.1 of \cite{Miller2}; for us it will be Lemma \ref{lem:globalminab} (which we use four times).
	
	The topic below is a favorite of the third named author, who has used it in several of his classes at Williams and at many math camps, including one attended by the other two. The problem is easy to state and explore, there are several natural generalizations, and when investigating numerous technical obstructions arise, making it a wonderful introduction to proofs and research. Similar to many problems in the field, what matters most is not an extensive background of advanced material but rather an inquisitiveness and ability to see connections and apply techniques. The Fibonacci numbers and other sequences arising from recurrence relations provide a rich ground for introducing students to research; after describing the problem, challenges and solutions, we end with some additional questions and encourage interested readers in pursuing these further to contact us.

	\section{Introduction}
	
	The \textit{egg drop problem} is a classic puzzle that has fascinated mathematicians and computer scientists. It is often featured in technical interviews at major technology companies such as Microsoft and Google, as it elegantly combines combinatorial reasoning with dynamic programming techniques.
	
	The Egg Drop problem was introduced in \cite{KonhauserVellemanWagon}.
	
	\begin{quote}
\textit{
	Suppose we wish to know which windows in a 36-story building are safe to drop an egg from, and which will cause the egg to break on landing. Suppose two eggs are available. What is the least number of egg droppings that is guaranteed to work in all cases?
}
\end{quote}

We make a few assumptions (see also Remark \ref{rek:worsecase}).
\begin{itemize}
\item An egg that survives a fall can be used again, with no damage.
\item A broken egg must be discarded.
\item The effect of a fall is the same for all eggs.
\item If an egg breaks when dropped, then it would break from a higher floor.
\item If an egg survives a fall, then it would survive a shorter fall.
\item It is not ruled out that the first-floor break eggs, nor is it ruled out that a drop from the 36th floor does not cause an egg to break.
\end{itemize}

The book displays the optimal solution to the problem -- to start at the eighth floor and if the egg breaks return to test every point from the first floor until the second egg breaks; otherwise use the first egg to test the fifteenth floor, and so on (as shown in Table \ref{intro}). No matter where the first egg breaks, it will always require no more than \(8\) drops.

\begin{table}[h!]\label{intro}
\centering
\begin{tabular}{c l}
	\toprule
	First Egg Breaks at Floor & Second Egg Test Floors if First Egg Breaks\\
	\midrule
	8 &  1 $\rightarrow$ 2 $\rightarrow$ 3 $\rightarrow$ 4 $\rightarrow$ 5 $\rightarrow$ 6 $\rightarrow$ 7 \\
	15 & 9 $\rightarrow$ 10 $\rightarrow$ 11 $\rightarrow$ 12 $\rightarrow$ 13 $\rightarrow$14 \\
	21 & 16 $\rightarrow$ 17 $\rightarrow$ 18 $\rightarrow$ 19 $\rightarrow$ 20 \\
	26 & 22 $\rightarrow$ 23 $\rightarrow$ 24 $\rightarrow$ 25 \\
	30 & 27 $\rightarrow$ 28 $\rightarrow$ 29 \\
	33 & 31 $\rightarrow$ 32 \\
	35 & 34 \\
	36 &  \\
	\bottomrule
\end{tabular}
\caption{Optimal dropping schedule for two eggs and 36 floors.}
\end{table}

Later, Boardman \cite{Boardman} generalized this into the problem: \textit{Given \(k\) eggs, how many floors can we handle if we have at most \(n\) drops?}

Using both generating functions and a direct counting approach, he showed that the total number of floors that can be tested is
\[
\sum_{j = 1}^{k} \binom{n}{j}.
\]

Building off of Boardman's work, Parks and Wills \cite{ParksWills} extended the problem using a binary decision tree to calculate the most floors that can be explored with \(k\) eggs and \(d\) drops for two adapted problems.

\vspace{1em}

\begin{itemize}
\item \textbf{Replacement Eggs}: the supply of eggs is restored to the original number \(k\) whenever the egg that is dropped does not break.
\item \textbf{Bonus Eggs}: A new egg is received whenever the egg that is dropped does not break.
\end{itemize}

\vspace{1em}

All of these results employ a \textit{dynamic search} where each drop decision depends on the previous outcomes. There is another family of methods known as \textit{statistical} or \textit{fixed-step} strategies that we use in this paper.

There are also many other adaptations that can be made.
This paper builds on earlier lectures presented by the third named author in his Problem Solving class, Math 331 (Lecture 9, \cite{Miller3}), at Williams College, as well as at many math camps.
We generalize the traditional egg drop problem into a problem with \(N\) floors and \(k\) eggs and bound the optimal step size and least number of drops required under such a strategy.
Beyond the one-dimensional problem, we also extend the traditional one-dimensional egg drop problem by introducing higher dimensional generalizations. Instead of searching a point along a vertical line of floors, we explore versions set in a \textit{two-dimensional plane} and \textit{three-dimensional space}, where the critical condition depends on multiple coordinates, and then create and solve a recursive equation to build up the general solution from the simple initial case. Moreover, we continue this line of generalization to search for a \textit{critical line} instead of a point in a 2D plane.

Our main theorems are the following.

\vspace{1em}

\textbf{Theorem \ref{thm:1D}: One Dimensional Critical Point Problem:}
\textit{For \(k \ge 1\) eggs and \(N\) floors, the minimum number of drops required in the worst-case scenario satisfies}
\[
P_1(k) \ \le \ \lceil k \cdot N^{1/k} \rceil.
\]

\vspace{1em}

\textbf{Theorem \ref{thm:2Dcriticalpoint}: Two Dimensional Critical Point Problem:}
\textit{For \(k \ge 2\) eggs and a rectangular region of dimensions \(M \times N\), the minimum number of drops required in the worst-case scenario satisfies}
\[
P_2(k) \ \le \ \lceil (k - 1) \cdot (M + N)^{1/(k - 1)} \rceil.
\]

\vspace{1em}

\textbf{Theorem \ref{thm:3D}: Three Dimensional Critical Point Problem:}
\textit{For \(k \ge 3\) eggs and a cubic space of dimensions \(L \times M \times N\), the minimum number of drops required in the worst-case scenario satisfies}
\[
P_3(k) \ \le \ \lceil (k - 2) \cdot (L + M + N)^{1/(k - 2)} \rceil.
\]

\vspace{1em}

\textbf{Theorem \ref{thm:2Dx+y=VM1}: Two Dimensional Critical Line \(x + y = V\), Method One:}
\textit{We perform a recursive jump search along the main diagonal of the rectangular region of dimensions \(M \times N\). For \(k \ge 1\) eggs, the minimum number of drops required in the worst-case scenario under \textbf{Method One} satisfies}
\[
L_2^{(1)}(k) \ \le \ \lceil k \cdot (M + N)^{1/k} \rceil.
\]

\vspace{1em}

\textbf{Theorem \ref{thm:2Dx+y=VM2}: Two Dimensional Critical Line \(x + y = V\), Method Two:}
\textit{We combine a diagonal search in the rectangular region of dimensions \(M \times N\) with a final verification of the last uncertain point using the reserved egg. For \(k \ge 2\) eggs, the minimum number of drops required in the worst-case scenario under \textbf{Method Two} satisfies}
\[
L_2^{(2)}(k) \ \le \ \lceil (k - 1) \cdot M^{1/(k - 1)} \rceil.
\]

\vspace{1em}

In future work, we plan to extend our study of the two-dimensional critical line \(x + y = V\) to the more general case \(\alpha x + \beta y = V\), which introduces new challenges and strategies.

\begin{remark}\label{rek:worsecase}
We focused on minimizing the worst-case number of drops, but we could also consider minimizing the average-case performance. Since the worst case occurs only rarely, it might be reasonable to allow some risk by optimizing the average instead. We do not address minimizing the average case in this paper, but leave it as a direction for future research.
\end{remark}


\section{Preliminaries}

Before moving to the one-dimensional case and introducing our recursive searching strategy, we state a general minimizing lemma that will be applied repeatedly in later sections. This lemma captures the common structure of the functions that appear subsequently in all \(1\)-, \(2\)- and \(3\)-dimensional cases when optimizing the step size, and hence the minimal number of drops.

\begin{lemma}\label{lem:globalminab}
Given any real numbers \(a > 0, b\ge 0, n > 1\), and any function of the form
\[
f(S) \ = \ \frac{n}{S} + a \cdot \left( \frac{n + b}{n} \right)^{1/a} \cdot S^{1/a}, \text{ for } S \in (0, n],
\]
the function always achieves its global minimum at
\[
S^* \ = \ n(n + b)^{-1/(a + 1)},
\]
while
\[
f(S^*) \ = \ (a + 1)(n + b)^{1/(a + 1)}.
\]

\end{lemma}

\begin{proof}
To minimize \(f(S)\), we calculate the first derivative
\begin{equation}\label{firstderivative}
	f'(S) \ = \ -n \cdot S^{-2} + \left( \frac{n + b}{n} \right)^{1/a} \cdot S^{(1 - a)/a}.
\end{equation}

Solving \(f'(S) = 0\), we find the critical point
\[
S^* \ = \ n(n + b)^{-1/(a + 1)}.
\]
(All the elementary calculations for the expressions derived in this section are provided in Appendix \ref{calculations}).
This is a potential location for a global minimum. Computing the second derivative at the critical point yields
\begin{equation}\label{secondderivative}
	\begin{aligned}
		f''(S^*) \ 
		&= \ 2n \cdot (S^*)^{-3} + \frac{1 - a}{a} \left( \frac{n + b}{n} \right)^{1/a} \cdot (S^*)^{(1 - 2a)/a}\\[6pt]
		&= \ n \cdot \frac{1 + a}{a} \cdot (S^*)^{-3} \ > \ 0,
	\end{aligned}
\end{equation}
since \(a, n, \text{ and } S^*\) are all positive. This verifies that the critical point is a unique local minimum. Since the domain is \(S \in (0, n]\), and \(S^*\) is not an endpoint, the global minimum must occur either at this interior point (i.e., the critical point) or at an endpoint.

\begin{itemize}
	\item As \( S \to 0^+ \):
	We have
	\(n/S \to +\infty \text{ and } a \cdot ((n + b)/n)^{1/a} \cdot S^{1/a} \to 0^+\). This implies that \(f(S) \to +\infty\).
	Thus, the function cannot achieve the global minimum there.
	\item At \(S = n\): We get
	
	\begin{equation}\label{f'(n)}
		f'(n) \ = \ -n \cdot n^{-2} + \left( \frac{n + b}{n} \right)^{1/a} \cdot n^{(1 - a)/a} \ = \ \frac{(n + b)^{1/a} - 1}{n} \ > \ 0,
	\end{equation}
	since the conditions \(a > 0, b \geq 0,\) and \( n > 1 \) imply that \((n + b)^{1/a} > 1\). Thus, the function cannot attain its minimum at the endpoint \(n\), as it increases at that point.
\end{itemize}
From the argument above, we can conclude that the critical point is a global minimum.

Substituting the optimal \(S^*\) into \(f(S)\), we obtain
\begin{equation}\label{f(S*)}
	\begin{aligned}
		f(S^*) \ 
		&= \ \frac{n}{S^*} + a \left( \frac{n + b}{n} \right)^{1/a}  (S^*)^{1/a} \\[6pt]
		&= \ (a + 1) \cdot (n + b)^{1/(a + 1)},
	\end{aligned}
\end{equation}
therefore completing the proof.
\end{proof}

\section{The One-Dimensional Critical Point Problem}

This section reviews a recursive strategy for the one-dimensional Egg Drop problem with \(k\) eggs and \(N\) floors. Distinct from the approach presented in the introduction, our method derives the optimal step size for the search and establishes an upper bound on the worst-case number of drops.\footnote{The worst-case number of drops ensures that, regardless of where the eggs break, the critical point can always be determined within this number of trials.} All results in this section are proved by induction, and our one-dimensional search also serves as the foundation for the higher-dimensional generalizations.

Among a linear array of numbers from \( 1 \) to \( N \), there exists a unique \emph{critical point} \( n \in (0, N] \) (both \(n\) and \(N\) are integers) such that

\begin{itemize}
\item if a magical egg is dropped from any point \( a < n \), the egg remains intact with no damage, and
\item if it is dropped from any point \( a \geq n \), the egg breaks.
\end{itemize}

We first describe a general strategy for solving this problem. If only one egg is available, test every lattice point from \( 1 \) upwards until the egg breaks at the critical point,\footnote{This method is known as a linear search or sequential search.} which is the only reliable method to determine the critical point.

\begin{figure}[h]
\begin{center}
	\scalebox{.21}{\includegraphics{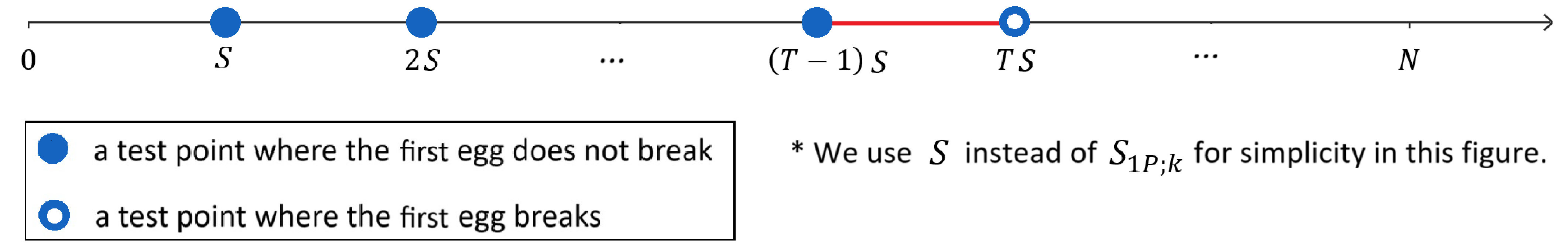}}
	\caption{\label{fig:Point1DGeneralCase} The 1D critical point strategy.}
	\end{center}\end{figure}

If more than one egg is available, we can introduce a more efficient searching approach. Suppose we have \( k \) eggs and \( N \) floors. We begin by dropping the first egg at positions of the form \( i \cdot S_{1P; k} \), where \( i = 1, 2, 3, \ldots \). Here \( S_{1P; k} \) denotes the \(k\)-egg step size to be determined when performing the \textit{critical point} (P) search in the one-dimensional setting.\footnote{This method  in known as a jump search.} Later we similarly use \( S_{2P; k} \) and \( S_{3P; k} \) for critical point searches in higher dimensions, and \( S_{2L; k} \) for the two-dimensional \textit{critical line} (L) search. Continue this process until the egg breaks at position \( T \cdot S_{1P; k} \). Now we know that the critical point must lie within the interval \( ((T - 1) S_{1P; k} , T	S_{1P; k}] \) (the red segment in Figure \ref{fig:Point1DGeneralCase}), consisting of exactly \( S_{1P; k} \) candidate points, among which \( S_{1P; k} - 1 \) need to be tested.\footnote{ Point \( T	S_{1P; k} \) does not need to be tested again.}

	We now face a reduced sub-problem: a \((k - 1)\)-egg search over a smaller interval of \( S_{1P; k} - 1 \) points. In the following recursive process, one egg is used to narrow down the interval, then the same strategy is applied with the remaining \( k - 1 \) eggs, until only one egg is left. At this point, we revert to a linear search within the final sub-segment.
	
	\begin{theorem} \label{thm:1D}
If there are \(k\) eggs, then in the worst-case scenario, the optimal step size under the given strategy is
\[
S_{1P; k} \ \le \ N^{(k - 1)/k},
\]
which leads to a minimum number of egg drops
\[
P_1(k) \ \le \ \lceil k \cdot N^{1/k} \rceil , \quad \text{for } k \geq 1.
\]
\end{theorem}

\subsection{Base Case: \( k = 1 \)}

In the worst-case scenario with only one egg, the egg does not break during the sequential search until it is dropped at position \(N\). Thus, the step size must be \(1\), and the number of drops required in the worst case is exactly \(N\).

The base case holds, since the optimal step size is \(S_{1P; 1} \le N^{(1 - 1)/1} = 1\), and the worst case number of drops is \(P_1(1) \le \lceil 1 \cdot N^{1/1} \rceil = N\).

\subsection{Inductive Case}

Assuming that Theorem \ref{thm:1D} is true for \(k\), we must show that \(S_{1P; k + 1} \le N^{k/(k + 1)}\) and \(P_1(k+1) \le \lceil (k + 1) \cdot N^{1/(k + 1)} \rceil \) are also true.

\begin{figure}[h]
\begin{center}
	\scalebox{.21}{\includegraphics{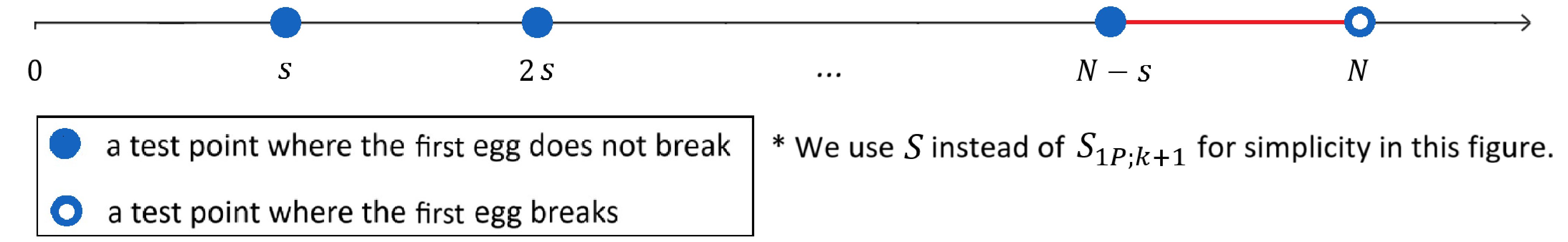}}
	\caption{\label{fig:Point1DWorstCase} The 1D critical point strategy in the worst case.}
	\end{center}\end{figure}

To do so, we apply the recursive strategy described previously. Suppose we drop the first egg at intervals of increment size \( S_{1P; k + 1} \) (i.e., jump search from positions \( S_{1P; k + 1}, 2S_{1P; k + 1}, \ldots \)), stopping when the egg breaks. In the worst case, this requires \( N / S_{1P; k + 1} \) drops with the first egg. Then, we are left with a sub-interval  (the red sub-segment in Figure \ref{fig:Point1DWorstCase}) of length at most \( S_{1P; k + 1} - 1 \) candidate points to search using the remaining \( k \) eggs.

	By the inductive hypothesis, the worst-case number of drops needed to search a segment of size \( S_{1P; k + 1} - 1 \) with \( k \) eggs is at most \( k \cdot (S_{1P; k + 1} - 1)^{1/k} \), which is approximately equal to, and less than, \( k \cdot S_{1P; k + 1} ^ {1/k} \) drops.
	
	Thus, the total number of drops in the worst case is approximately and at most
	\begin{equation} \label{eq:calculation}
f_{1P; k + 1}(S_{1P; k + 1}) \ \leq \ \frac{N}{S_{1P; k + 1}} + k \cdot S_{1P; k + 1} ^{1/k}.
\end{equation}

We now use Lemma \ref{lem:globalminab} with \(a = k,\ b = 0, \text{ and } n = N\) giving
\[
S_{1P; k + 1}^* \ = \ N \cdot (N + 0)^{-1/(k + 1)} \ = \ N^{k/(k + 1)}
\]
and
\[
f_{1P; k + 1}(S_{1P; k + 1}^*) \ = \ (k + 1)(N + 0)^{1/(k + 1)} \ = \ (k + 1) \cdot N^{1/(k + 1)}.
\]

Note that in Lemma \ref{lem:globalminab} to minimize \(f(S)\), we treated \(S\) as a continuous variable to allow differentiation. However, in a real-world setting, such as a building, the tested points correspond to discrete floors and must therefore be integers, so we round up to the nearest integer floor. In an abstract setting, such as a number line, non-integer points can be tested.

In this case, the number of drops is taken as \(\lceil (k + 1) \cdot N^{1/(k + 1)} \rceil\) to obtain an integer value. The approximation is reasonable and safe because the derivative of \(f_{1P; k+1}(S_{1P; k + 1}) \) decreases and then increases, so we know that the best integer step size lies either immediately above or below the continuous optimal. And there is no problem with rounding UP the step size to the closest integer, since there is only a very small change in the total number of required drops between rounding up and down.\footnote{For example, with \(N = 100\) and \(k = 4\), the optimal continuous step size is approximately \(S_{1P; 4} \approx 31.62278\), giving \(f_{1P; 4}(31.62278) \approx 12.34778\) drops. Rounding up to \( 32\) yields \(f_{1P; 4}(32) \approx 12.63866\), while rounding down to \(31\) gives \(f_{1P; 4}(31) \approx 12.66425\). Since we're only aiming for an upper bound, choosing the ``wrong'' integer step size merely increases the bound by a negligible amount, and both round to the same integer.} However, as illustrated by the \textbf{Knapsack Problem} (see, for example, \cite{Fr, Miller1}), the integer minimum or maximum can significantly differ from its real-valued counterpart.\footnote{Imagine a bag that can hold only 100 pounds, and we want to maximize the value of what is inside; assume we can put any amounts of the following three products in it, so long as the total weight is at most 100. One product is worth \$150 and comes in 51 pound units, the second is \$100 per unit and each is 50 pounds, while the last is worth \$99 and comes in 50 pound units. Clearly we never should use any of the third item. If we can take partial amounts of a product (i.e., we are in a continuous case), we take 100/51 units of the first product and none of the second, for a value of about \$292.12. If however we can only take integral amounts of each product, the best we can do is to take none of the first and two of the second, for a value of \$200. Note the integral optimal is far from the real.} So, we need to be cautious when rounding.

Thus, given \( S_{1P; k} \le N^{k - 1/k} \) is true, \( S_{1P; k + 1} \le N^{k/(k+1)}\) is true, and given \( P_1(k) \le \lceil k \cdot N^{1/k} \rceil \) is true, \( P_1(k+1) \le \lceil (k + 1) \cdot N^{1/(k + 1)} \rceil \) is true, completing the inductive step, and hence the proof.
\hfill $\Box$

\vspace{1em}

Since we approximated \(S_{1P; k + 1} - 1\) into \(S_{1P; k + 1}\), and rounded up the minimum number of drops, we obtain our claimed upper bound.

Also, note that this is not the optimal testing strategy, though it is close. When \(k = 2\) and the total number of floors is a triangular number such as \(36\), the method described in the introduction achieves a slightly better result. In that approach, the search steps decrease by one each time, ensuring that the total number of drops remains the same, no matter where the first egg breaks. However, in our theorem, the step sizes are uniform, resulting in the same power of \(N\) as we get with the "triangular" strategy, which makes the analysis easier and the generalization to higher-dimensional settings more natural.

\section{The Two-Dimensional Critical Point Problem}

There are many ways to extend the traditional egg drop problem, and generalizing into a two-dimensional critical point problem is a natural one. Let us consider a rectangular region in the plane with coordinates ranging from \( (0, 0) \) to \( (M, N) \), without loss of generality assuming \( M \geq N \). There exists a unique \emph{critical point} \( (m, n) \), where \(m\), \(n\), \(M\), and \(N\) are all integers, with \(m \in (0, M]\) and \(n \in (0, N]\) such that if a magical egg is dropped at a point \( (a, b) \), then

\begin{itemize}
\item if \( a < m \) and \( b < n \), the egg remains intact, and undamaged;
\item otherwise (i.e., if \( a \geq m \) or \( b \geq n \)), the egg breaks.
\end{itemize}

The objective is to determine the exact value of the critical point \( (m, n) \) using the minimal number of egg drops.

It is impossible to determine both coordinates of the critical point using only one egg, as a single egg can only be used to resolve the critical coordinate along one axis.

If two eggs are available, we can
only find the critical point by performing a linear search with one egg per axis, until the egg breaks along each axis.

\begin{figure}[h]
\begin{center}
	\scalebox{.09}{\includegraphics{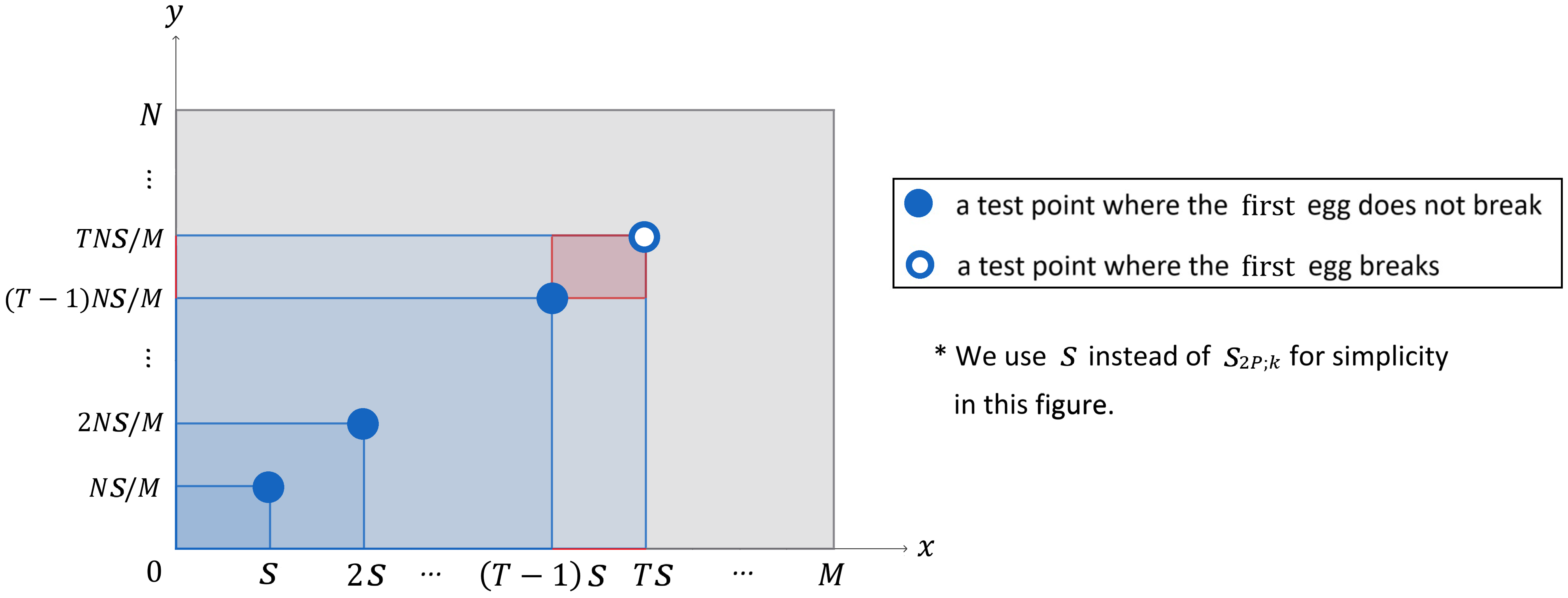}}
	\caption{\label{fig:Point2DGeneralCase} The 2D critical point strategy.}
\end{center}
\end{figure}

When more than two eggs are provided, a more efficient strategy may be employed. Suppose we have \( k \geq 2 \) eggs. We begin our jump search by dropping the first egg at positions of the form \( (iS_{2P; k}, iNS_{2P; k}/M) \), where \( i = 1, 2, 3, \ldots \) and \( S_{2P; k} \) is a \(k\)-egg step size (or a scaling factor) to be determined in the two-dimensional setting when searching for a \textit{critical point} (P) (as shown in Figure \ref{fig:Point2DGeneralCase}).\footnote{The rectangles defined between \((0,0)\) and \((iS_{2P; k}, iNS_{2P; k}/M)\) are all similar to the largest rectangular region, and all the test points lie along the diagonal of this rectangle.} Note that in the real-world context, we would round \( iNS_{2P; k}/M \) to the nearest integer to ensure that drops occur at integer positions. However, in an abstract setting, such as a coordinate plane, we allow the egg to be dropped at any real-valued position, so \( iNS_{2P; k}/M \) need not be an integer. Continue this process until the egg breaks at position \( (TS_{2P; k}, TNS_{2P; k}/M) \). Now, the critical point must lie within the red sub-rectangle from \( ((T - 1)S_{2P; k}, (T-1)NS_{2P; k}/M)\) to \( (TS_{2P; k}, TNS_{2P; k}/M) \). This red rectangle has side lengths \( S_{2P; k} \) and \( NS_{2P; k}/M \), respectively.

Having lost one egg, this becomes a reduced sub-problem: a \((k - 1)\)-egg search within the red sub-rectangle of dimensions \( S_{2P; k} \) by \( NS_{2P; k}/M \). Continue the recursive approach until only two eggs are left. Then, we fall back on the two-egg strategy described earlier to find the critical point.

\begin{theorem}\label{thm:2Dcriticalpoint}
If there are \(k\) eggs, then, under the strategy described above, in the worst-case scenario, the minimum number of drops is given by
\[
P_2(k) \ \leq \ \lceil (k - 1) \cdot (M + N)^{1/(k - 1)} \rceil, \quad \text{for } k \geq 2.
\]
\end{theorem}

\subsection{Base Case: \( k = 2 \)}

When we have two eggs, the only guaranteed strategy is to perform a linear search along each axis---one egg used to determine the critical \( x \)-coordinate and the other for \( y \). In the worst-case scenario, the eggs do not break until it is dropped at position \( M \) along the \( x \)-axis and at position \( N \) along the \( y \)-axis. This results in a total of \( M + N \) egg drops.

The base case is true, since the worst case number of drops is \(P_2(1) \le \lceil (2 -1) \cdot (M + N)^{1/(2 - 1)} \rceil = M + N \).

\subsection{Inductive Step}

Assuming that Theorem \ref{thm:2Dcriticalpoint} is true for \(k\), then we must show that \(P_2(k + 1) \le \lceil k \cdot (M + N)^{1/k} \rceil \) is also true.

\begin{figure}[h]
\begin{center}
	\scalebox{.14}{\includegraphics{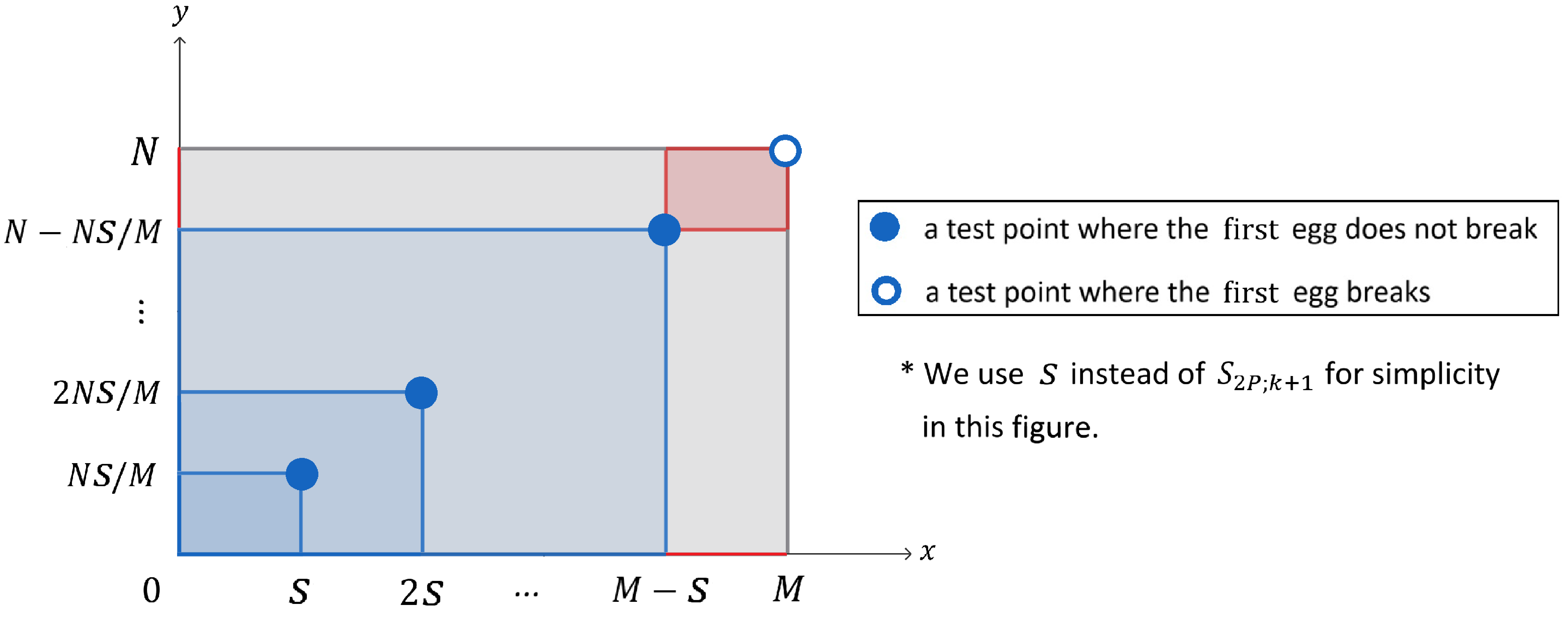}}
	\caption{\label{fig:Point2DWorstCase} The 2D critical point strategy in the worst case.}
	\end{center}\end{figure}

To do this, we consider the same strategy described earlier: we drop the first egg at points of the form \( (iS_{2P; k + 1}, i NS_{2P; k + 1}/M) \), for \( i = 1, 2, 3, \ldots \), until the egg breaks, to narrow down the searching region (as shown in Figure \ref{fig:Point2DWorstCase}). Again, \( iNS_{2P; k + 1}/M\) can be any real number.  In the worst case, this requires \( M / S_{2P; k + 1} \) drops with the first egg.

Once the egg breaks, we are left with a reduced search region: the red sub-rectangle, with a diagonal from \( (M - S_{2P; k + 1}, N - NS_{2P; k + 1}/M) \) to \( (M, N) \), and a total side length of
\[
S_{2P; k + 1} + \frac{NS_{2P; k + 1}}{M} \ = \ \frac{M + N}{M} \cdot S_{2P; k + 1}.
\]

	Applying the inductive hypothesis to this smaller region with the remaining \( k \) eggs, the worst-case number of drops is
	\[
	(k - 1) \cdot \left( \frac{M + N}{M} \cdot S_{2P; k + 1} \right) ^{1/(k - 1)} \ = \ (k - 1) \cdot \left(\frac{M + N}{M} \right)^{1/(k - 1)} \cdot S_{2P; k + 1}^{1/(k - 1)}.
	\]
	
	Thus, the total number of drops is the sum of the diagonal jump search and the linear search, and can be expressed as function
	\[
	f_{2P; k + 1}(S_{2P; k + 1}) \ \le \ \frac{M}{S_{2P; k + 1}} + (k - 1) \cdot \left( \frac{M + N}{M} \right)^{1/(k - 1)} \cdot S_{2P; k + 1}^{1/(k - 1)}.
	\]
	
	We now use Lemma \ref{lem:globalminab} as \(a = k - 1, b = N, \text{ and } n = M\) giving,
	\[
	S_{2P; k + 1}^* \ = \ M \cdot (M + N)^{-1/k},
	\]
	and
	\[
	f_{2P; k + 1}(S_{2P; k + 1}^*) \ = \ k \cdot (M + N)^{1/k}.
	\]
	
	Therefore, it is true that the minimum number of egg drops required in the worst case with \( k + 1 \) eggs is
	\[
	P_2(k + 1) \ \le \ \lceil k \cdot (M + N)^{1/k} \rceil,
	\] when given \(P_2(k) \le \lceil (k - 1) \cdot (M + N)^{1/(k - 1)} \rceil \) is true, completing the inductive step and thus the proof.
	\hfill $\Box$
	
	\vspace{1em}
	
	The generalization of the critical point to the three-dimensional case is similar to the two-dimensional case. For the full argument, see Appendix \ref{appendix:3D-critical}.
	
	The results in one, two, and three dimensions suggest a clear recursive pattern. This observation naturally leads to the following conjecture for a general \(d\)-dimensional setting.
	
	\begin{conjecture}
In a \(d\)-dimensional setting with sides \(N_1, N_2, \dots, N_d\) and \(k\) eggs, the minimum number of drops required in the worst-case scenario under the same strategy satisfies
\[
P_d(k) \ \leq \ \lceil (k - d + 1) \cdot (N_1 + N_2 + \dots + N_d)^{1/(k - d + 1)} \rceil, \quad \text{for } k \geq d.
\]
\end{conjecture}

The tricky part in this expression is figuring out the correct coefficient and the denominator of the exponent. In the one-dimensional setting, both are \(k\), or \(k - 0\), and in the two-dimensional case, they become \(k - 1\). We see a pattern, and thus, in the \(d\)-dimensional case, we use \(k - (d - 1)\) (as we subtract one less than the dimension).

\section{The Two-Dimensional Critical Line \(x + y = V\) Problem} \label{critical line}

This section introduces a new egg drop problem and a different rule for classifying whether eggs break or not when dropped.
Consider a rectangular region in the plane with coordinates ranging from \( (0, 0) \) to \( (M, N) \). Without loss of generality, assume \( M \geq N \). Within this region, there exists a unique \emph{critical line} of the form \( x + y = V \) that passes through the rectangle, where \(x\) and \(y\) are variables, and \(V\) is an unknown positive constant. The critical line determines the outcome of dropping a magical egg at a point \( (a, b) \):

\begin{itemize}
\item if the point \((a, b)\) lies below the line (i.e., \( a + b < V\)), the egg remains intact with no damage;
\item otherwise (i.e., if \( a + b \geq V\)), the egg breaks.
\end{itemize}

The objective is to split the \((M + 1)(N + 1)\) lattice points within and on the side of the rectangle into two groups: safe points and breaking points.\footnote{Safe points are those where the egg remains intact, while breaking points are those where the egg breaks.} To determine the location of a line we usually need two points the line passes through. However, since we know the slope of the line is \(-1\), we only need to locate one point the line passes through.

\begin{figure}[h]
\begin{center}
	\scalebox{.08}{\includegraphics{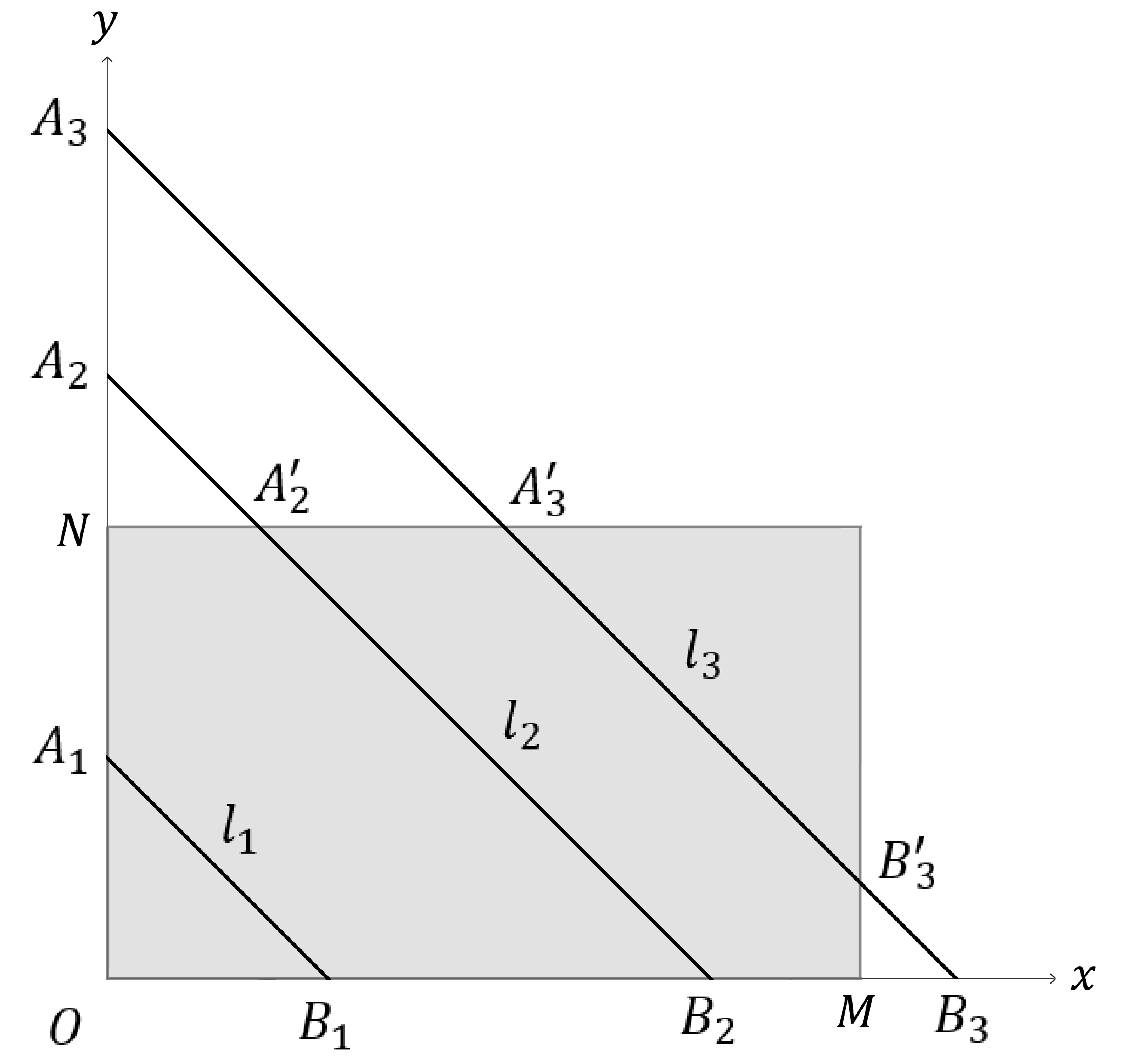}}
	\caption{\label{fig:LineXY1EggSame} The 2D critical line \(x + y = V\): 1-egg argument.}
	\end{center}\end{figure}

Let's start by outlining our overall strategy. With only one egg, we need to take the most cautious approach---a linear search along the bottom and right edges of the rectangle---to guarantee that we do not skip over any lattice points on the critical line, since we have no spare eggs to go back and test again. It is the same whether we move along the bottom and right sides, or along the left and top sides, since the slope of the critical line is \(-1\) (as shown in Figure \ref{fig:LineXY1EggSame}).

\begin{itemize}
\item If the critical line \(l_1\) intersects the left and bottom sides of the rectangle at \(A_1\) and \(B_1\), respectively, then \(OA_1 = OB_1\).

\item If the critical line \(l_2\) intersects the top and bottom sides at \(A'_2\) and \(B_2\), and the \(y\)-axis at \(A_2\), then \(NA_2 = NA'_2\) and \(OA_2 = OB_2\). This implies that \(ON + NA'_2 = OB_2\).

\item If the critical line \(l_3\) intersects the top and right sides at \(A'_3\) and \(B'_3\), and the \(y\) and \(x\)-axes at \(A_3\) and \(B_3\), respectively, then \(NA_3 = NA'_3\), \(MB_3 = MB'_3\), and \(OA_3 = OB_3\). This implies that \(ON + NA'_3 = OM + MB'_3\).

\end{itemize}

	With multiple eggs, more efficient strategies can be implemented.

	\subsection{Method One}
	
\begin{figure}[h]
\begin{center}
	\scalebox{0.08}{\includegraphics{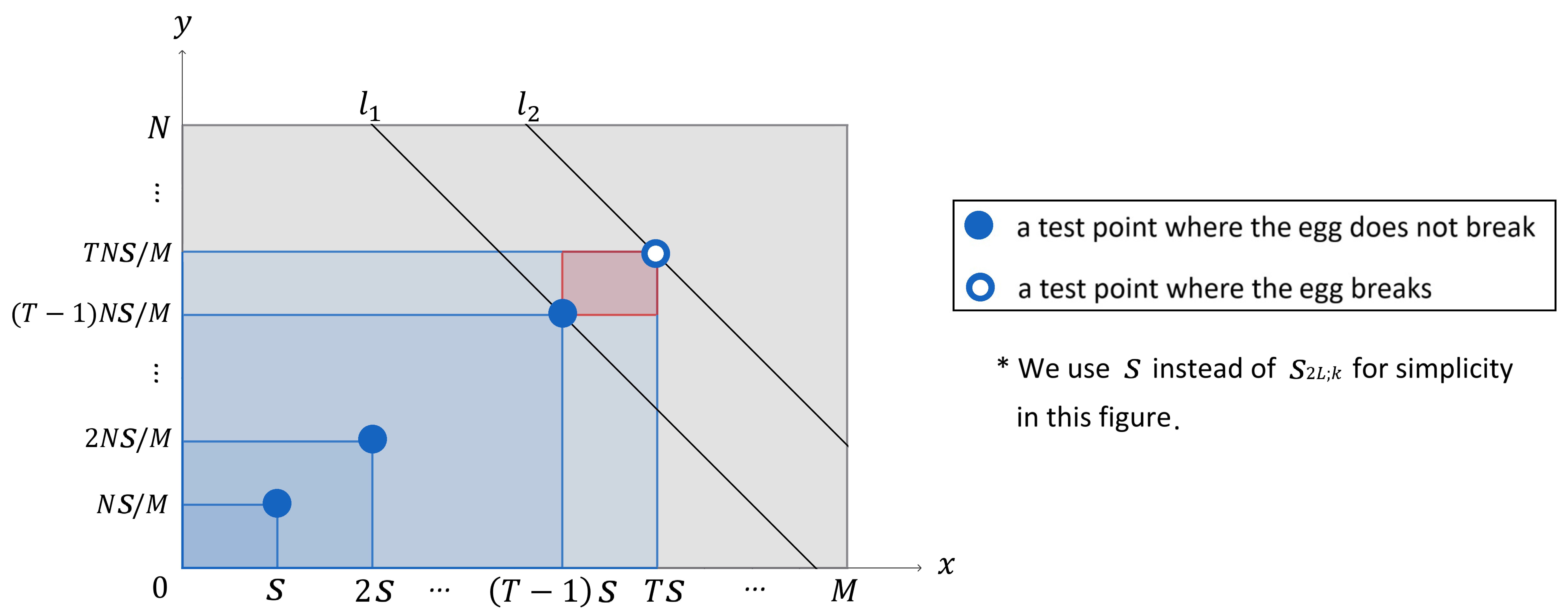}}
	\caption{\label{fig:LineXYkEgg} The 2D critical line \(x + y = V\): Method One strategy.}
	\end{center}\end{figure}

	Suppose we have \( k \geq 2 \) eggs. We use Figure \ref{fig:LineXYkEgg} to illustrate our first strategy here. We begin our jump search by dropping the first egg at points of the form \( (iS_{2L; k}, iNS_{2L; k}/M) \), where \( i = 1, 2, 3, \dots \) along the diagonal of the \(M\)-by-\(N\) rectangle, and \( S_{2L; k} \) is a scale factor to be determined when searching for the \textit{critical line} (L) in the 2D setting. Assuming the same convention as above, \(iNS_{2L; k}/M\) can be any real numbers in the abstract setting. Continue this process until the egg breaks at point \( (TS_{2L; k}, TNS_{2L; k}/M) \). Now, the critical line must pass through the red sub-rectangle from \( ((T - 1)S_{2L; k}, (T-1)NS_{2L; k}/M)\) (excluded) to \( (TS_{2L; k}, TNS_{2L; k}/M) \). We know that this red rectangle has side lengths \( S_{2L; k} \) and \( NS_{2L; k}/M \), respectively, and
	
	\begin{itemize}
\item all points on or to the left of the line \(l_1\) are safe points, and
\item all points on or to the right of the line \(l_2\) are breaking points.
\end{itemize}

	Having lost one egg, we are reduced to a sub-problem: a \((k - 1)\)-egg search within the red sub-rectangle of dimensions \( s \) by \( Ns/M \). Continue the recursive approach along the diagonal until only one egg remains. Then, we fall back onto the one-egg strategy described earlier.
	
	\begin{theorem} \label{thm:2Dx+y=VM1}
If there are \(k\) eggs, then, under Method One, in the worst-case scenario, the minimum number of drops is given by
\[
L_2^{(1)}(k) \ \leq \ \lceil k \cdot (M + N)^{1/k} \rceil, \quad \text{for } k \geq 1.
\]
\end{theorem}

\subsubsection{Base Case: \(k = 1\)}

\begin{figure}[h]
\centering
\scalebox{.15}{\includegraphics{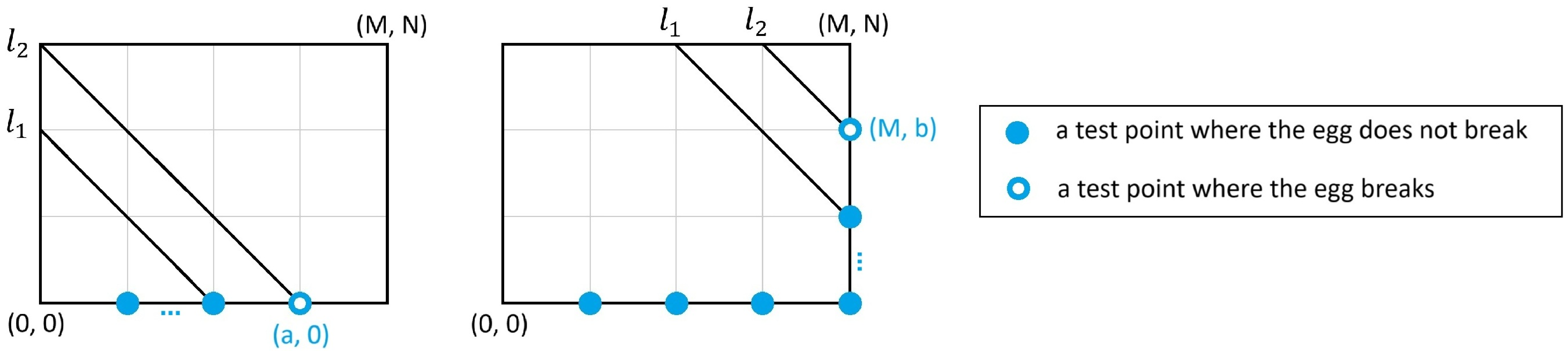}}
\caption{\label{fig:LineXY1Egg} The 2D critical line \(x + y = V\): 1-egg strategy.}
\end{figure}

If the egg breaks at point \((a,0)\) on the bottom side of the rectangle (see Figure \ref{fig:LineXY1Egg}),

\begin{itemize}
\item all points on or to the left of the line \(l_1: y = -x + a - 1\) are safe points, and
\item all points on or to the right of the line \(l_2: y = -x + a\) are breaking points.
\end{itemize}

The line \(l_2\) is simply \(l_1\) shifted one square to the right. Since \(l_1\) has slope \(-1\) and an integer \(x\)-intercept, there are no lattice points between \(l_1\) and \(l_2\).

Similarly, if the egg breaks at point \((M,b)\) on the right side of the rectangle, there are no lattice points between \(l_1\) and \(l_2\).

In the worst case, with only one egg, we test all points on the bottom and right sides of the rectangle until the egg breaks at point \((M,N)\), with the total number of drops being \(M + N\).
Note that \(L_2^{(1)}(1) = \lceil 1 \cdot (M + N)^{1/1} \rceil = M + N\).
Thus, the base case holds at \(1\).

\subsubsection{Inductive Step}

Assuming that Theorem \ref{thm:2Dx+y=VM1} is true for \(k\), then we must show that \(L_2^{(1)}(k + 1) \le \lceil (k + 1) \cdot (M + N)^{1/(k + 1)} \rceil \) is true.

To do so, we consider the same strategy described earlier: we drop the first egg at points of the form \( (iS_{2L; k + 1}, i NS_{2L; k + 1}/M) \), for \( i = 1, 2, 3, \ldots \), until the egg breaks, to narrow down the searching region. Again, \(iNS_{2L; k + 1}/M\) can be any real numbers, not necessarily integers. In the worst case, this requires \( M / S_{2L; k + 1} \) drops with the first egg breaking at the very end of the diagonal of the rectangle.

Once the egg breaks, we are left with a reduced search region: the sub-rectangle, with a diagonal from \( (M - S_{2L; k + 1}, N - NS_{2L; k + 1}/M) \) to \( (M, N) \), and a total side length of
\[
S_{2L; k + 1} + \frac{NS_{2L; k + 1}}{M} \ = \ \frac{M + N}{M} \cdot S_{2L; k + 1}.
\]

Applying the inductive hypothesis to this smaller region with the remaining \( k \) eggs, the worst-case number of drops is
\[
k \cdot \left( \frac{M + N}{M} \cdot S_{2L; k + 1} \right) ^{1/k} \ = \ k \cdot \left(\frac{M + N}{M} \right)^{1/k} \cdot S_{2L; k + 1}^{1/k}.
\]

Thus, the total number of drops is the sum of the diagonal jump search and the linear search along the sides of the last sub-rectangle, and can be expressed as function
\[
f_{2L; k + 1}(S_{2L; k + 1}) \ \le \ \frac{M}{S_{2L; k + 1}} + k \cdot \left(\frac{M + N}{M} \right)^{1/k} \cdot S_{2L; k + 1}^{1/k}.
\]

We now use Lemma \ref{lem:globalminab} as \(a = k, b = N, \text{ and } n = M\) giving,
\[
S_{2L; k + 1}^* \ = \ M \cdot (M + N)^{-1/(k + 1)},
\]
and
\[
f_{2L; k + 1}(S_{2L; k + 1}^*) \ = \ (k + 1)(M + N)^{1/(k + 1)}.
\]

Therefore, when given \(L_2^{(1)}(k) \le \lceil k \cdot (M + N)^{1/k} \rceil \) is true, it is also true that the minimum number of egg drops required in the worst case with \( k + 1 \) eggs is
\[
L_2^{(1)}(k + 1) \ \le \ \lceil (k + 1) \cdot (M + N)^{1/(k + 1)} \rceil,
\] completing the inductive step. \hfill $\Box$

\subsection{Method Two}

With two eggs, we can also apply an alternative strategy: if we test points along the diagonal of the form \((i, iN/M)\), where \(i = 1,2,3, \dots\), with \(iN/M\) can be any real numbers, until the first egg breaks at point \((T, TN/M)\), then

\begin{itemize}
\item all points on or to the left of the line \(l_1: y = -x + (T - 1)(1 + N/M)\) are safe points, and

\item all points on or to the right of the line \(l_2: y = -x + T(1 + N/M)\) are breaking points.
\end{itemize}

Note that \(l_2\) is simply \(l_1\) shifted \(1 + N/M\) upwards.
So there can be one or two rows of points between \(l_1\) and \(l_2\), since \(M \geq N\) implies that \(0 \leq N/M \leq 1\) and \(1 \leq 1 + N/M \leq 2\).

\subsubsection*{\underline{Worst Case One: First Egg Breaks at \((M - 1,  (M - 1)N/M)\)}}

\begin{figure}[h]
\centering
\scalebox{.15}{\includegraphics{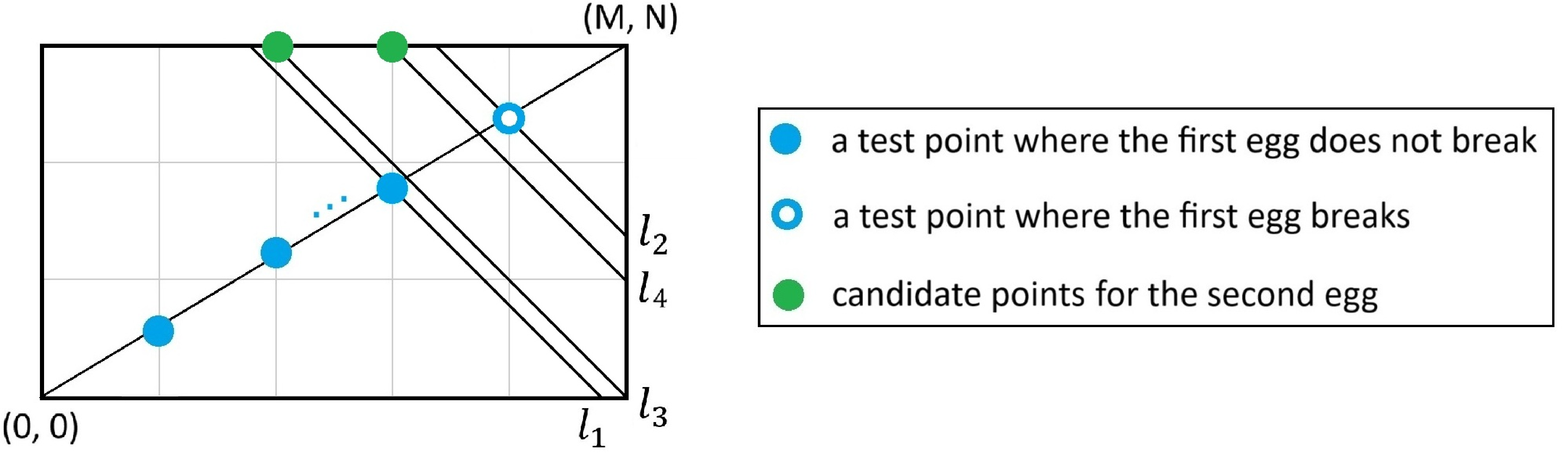}}
\caption{\label{fig:LineXY2EggWorstCase1} The 2D critical line \(x + y = V\): Method Two: 2-egg strategy in Worst Case One.}
\end{figure}

If the first egg breaks at point \((M - 1,  (M - 1)N/M)\) (as shown in Figure \ref{fig:LineXY2EggWorstCase1}), then it takes us \(M - 1\) drops, leaving two rows of candidate points between \(l_1\) and \(l_2\) on lines \(l_3: y = -x + M + N - 3\) and \(l_4: y = -x + M + N - 2\) in the worst case. We only need to test one point on each of \(l_3\) and \(l_4\) respectively with the second egg. Thus, altogether this process requires \((M-1) + 2 = M+1\)  drops.

\subsubsection*{\underline{Worst Case Two: First Egg Breaks at \((M, N)\)}}

\begin{figure}[h]
\centering
\scalebox{.15}{\includegraphics{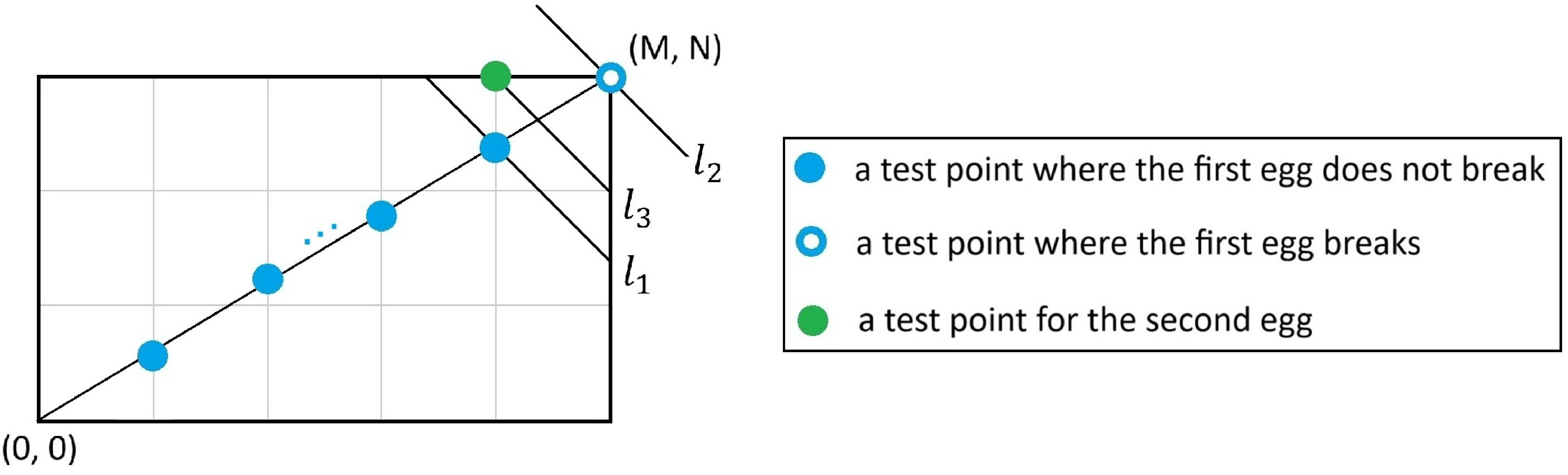}}
\caption{\label{fig:LineXY2EggWorstCase2} The 2D critical line \(x + y = V\): Method Two: 2-egg strategy in Worst Case Two.}
\end{figure}

If the first egg breaks at point \((M, N)\) (as shown in Figure \ref{fig:LineXY2EggWorstCase2}), it requires \(M\) drops, and leaves only one row of points between \(l_1\) and \(l_2\) on line \(l_3: y = -x + M + N - 1\). We only need to test one of these points on \(l_3\) with the second egg, at either point \((M - 1, N)\) or \((M, N - 1)\). Thus, altogether this process also requires \(M+1\)  drops.

\begin{theorem} \label{thm:2Dx+y=VM2}
Under Method Two, with two eggs available, the worst case occurs when the first egg breaks at \((M - 1,  (M - 1)N/M)\) or \((M,N)\), requiring \(M+1\)  drops.
\end{theorem}

With \(k \geq 3\) eggs, we first set one egg aside and use the remaining \(k - 1\) eggs to search along the diagonal of the rectangle. This requires \((k - 1) \cdot M^{1/(k - 1)}\) drops, since the diagonal can be treated as a line where we apply a one-dimensional searching strategy (Theorem \ref{thm:1D}). We then use the reserved egg for testing the final uncertain point, costing another drop in the worst case.

\begin{theorem}
If there are \(k\) eggs, then, under Method Two, in the worst-case scenario, the minimum number of drops is given by
\[
L_2^{(2)}(k) \ \leq \ \lceil (k - 1) \cdot M^{1/(k - 1)} \rceil + 1, \quad \text{for } k \geq 2.
\]
\end{theorem}

\section{Future Work}

Consider a rectangular region in the plane with coordinates ranging from \((0, 0)\) to \((M, N)\). Without loss of generality, assume \(M \geq N\). Within this region, there exists a unique critical line of the form \(\alpha x + \beta y = V\) that passes through the rectangle, where \(x\) and \(y\) are variables, \(\alpha\), \(\beta\) and \(V\) are unknown constants with following constraints.

\begin{itemize}
\item The critical line has a positive \(y\)-intercept.
\item The slope of the critical line can be negative, \(0\), or undefined, so that the line can be descending, horizontal or vertical.
\item The critical line has to go through at least two lattice points within or on the boundary of the rectangle.

\end{itemize}

The critical line determines the outcome of dropping a magical egg at a point \((a, b)\):
\begin{itemize}
\item if the point \((a, b)\) lies strictly below the line (i.e., \(\alpha x + \beta y < V\) ), the egg remains intact and undamaged;
\item otherwise (i.e., if \(\alpha x + \beta y \geq V\) ), the egg breaks.
\end{itemize}

The objective is to split the \((M + 1)(N + 1)\) lattice points on and within the rectangle into two groups: safe points and breaking points.

\subsection{Failure of One-Egg Strategy}

It is impossible to complete the classification with only one egg, because we need at least two breaking points to narrow down the range and constrain the slope of the critical line.

\begin{figure}[h]
\centering
\scalebox{.18}{\includegraphics{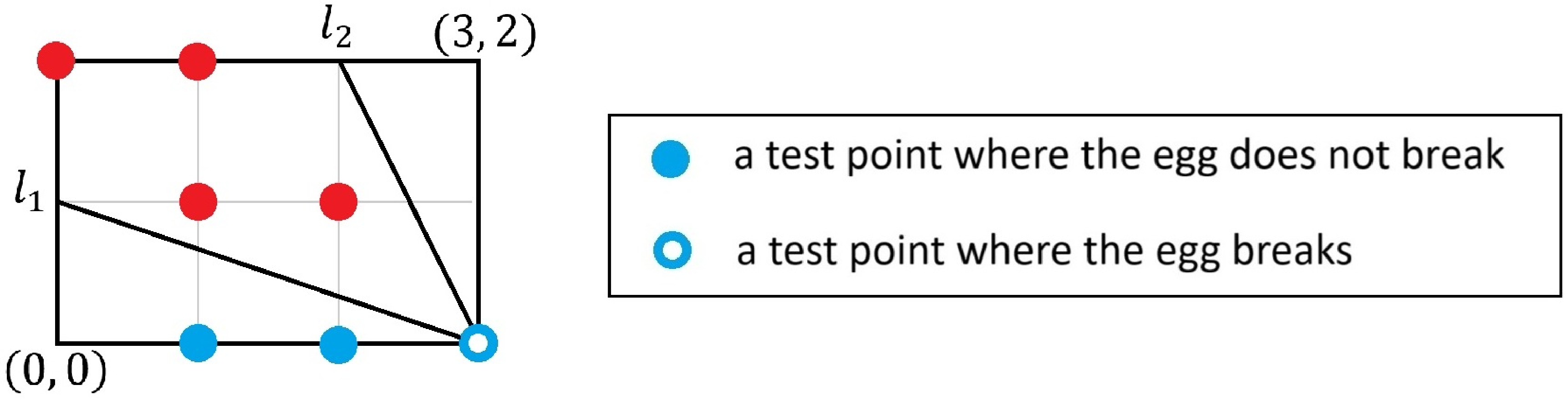}}
\caption{\label{fig:LineaXbY1Egg} The 2D critical line \(\alpha x + \beta y = V\): 1-egg failure example.}
\end{figure}

For example, in Figure \ref{fig:LineaXbY1Egg}, the egg breaks at point (3, 0). There are six possible locations for the critical line, two of which are shown as \(l_1\) and \(l_2\):
\begin{itemize}
\item \(l_1\) is a line that enters through point \((0, 1)\) and exits through point \((3, 0)\), or
\item \(l_2\) is a line that enters through point \((2, 2)\) and exits through point \((3, 0)\).
\end{itemize}

If \(l_1\) is the critical line, the set of safe points includes \((0, 0)\), \((1, 0)\) and \((2, 0)\), and the rest are breaking points. If \(l_2\) is the critical line, there will be four more safe points. Thus, with only one egg, we cannot determine the complete set of breaking and non-breaking points.

\subsection{Two-Egg Strategy}\label{SlopeUnknown2Egg}

With only two eggs, we must use the most conservative approach to ensure that no lattice points on the critical line are skipped.
\begin{itemize}
\item First, perform a sequential search at every lattice point along the left and top sides of the rectangle using the first egg, and suppose the egg breaks at point \((0,a)\) on the left side or point \((b, N)\) on top.
\item If the egg breaks at \((0, a)\), we then consider all points \((x,y)\) within or on the boundary of the rectangle from \((1, 0)\) to \((M, a)\) and compute their slopes relative to \((0,a)\), namely
\[
m_{1}(x,y) \ = \ \frac{y-a}{x}.
\]
If the egg breaks at \((b, N)\), we can calculate the slopes from \((b, N)\) to every point \((x, y)\) in the rectangle from \((b + 1, 0)\) to \((M, N)\), giving
\[
m_{2}(x,y) \ = \ \frac{y - N}{x - b}.
\]
\item The set of slopes is ordered in increasing value, and the corresponding points are tested sequentially from the smallest slope upward, using the second egg, until it breaks at a certain point \((c, d)\).
\item Because the critical line must pass through at least two lattice points, points that share the same slope respect to the first breaking point (either \((0,a)\) or \((b, N)\)) can be treated as a single candidate point.\footnote{Without the constraint that the critical line passes through at least two lattice points, even points with same slope relative to the first breaking point must be strictly tested from the lowest upward.}
\end{itemize}

At this point, we can conclude that all points strictly below the line passing through the first breaking point \((0,a)\) (or \((b,N)\)) and the second breaking point \((c,d)\) are safe, while all points on or above this line are points where the egg is bound to break.

From here we leave this critical line (unknown slope) problem for future work. In particular, we have not yet generalized the strategy to handle more than two eggs---such as cases of three, four, or \(k\) eggs---and have not derived a general formula to compute or bound the number of drops required. It also remains an open question whether alternative strategies could be employed, for instance, using smaller increments with more frequent drops when additional eggs are available.

\subsection{Other Questions}\label{subsec:otherquestions}

Another direction for future research is to consider average-case strategies. In this paper, we focus on minimizing the number of drops in the worst-case scenario, but we could also minimize the expected number of drops. Since the worst case occurs only rarely, it may be reasonable to allow some risk to improve average performance. If we assume a probability distribution for the breaking floor---for example, uniform, exponential, or binomial---then the goal would be to optimize the search strategy with respect to that distribution.

Additionally, we could compare the two methods and see when one is better than the other. We have done some preliminary calculations along these lines; see Appendix \ref{analysis}.


\appendix

\section{Elementary Calculations}\label{calculations}

We begin with the first derivative from Equation \ref{firstderivative}
\[
f'(S) \ = \ -n \cdot S^{-2} + \left( \frac{n + b}{n} \right)^{1/a} \cdot S^{(1 - a)/a}.
\]
We first set \(f'(S) = 0\) and get
\[
-n \cdot (S^*)^{-2} + \left( \frac{n + b}{n} \right)^{1/a} \cdot (S^*)^{(1 - a)/a} \ = \ 0.
\]
Add \(n \cdot (S^*)^{-2}\) to both sides and we have
\[
\left( \frac{n + b}{n} \right)^{1/a} \cdot (S^*)^{(1 - a)/a} \ = \ n \cdot (S^*)^{-2}.
\]
Multiplying both sides by \((S^*)^{2}\) gives us
\[
\left( \frac{n + b}{n} \right)^{1/a} \cdot (S^*)^{(a + 1)/a} \ = \ n.
\]
Multiply both sides by \((n/(n + b))^{1/a}\), and we get
\[
(S^*)^{(a + 1)/a} \ = \ n \cdot \left( \frac{n}{n + b} \right)^{1/a} \ = \ n^{(a + 1)/a} \cdot (n + b)^{-1/a}.
\]
Raise both sides to the power of \(a/(a + 1)\), yielding
\[
S^{*} = n(n + b)^{-1/(a + 1)}.
\]

Then we examine the second derivative from Equation \ref{secondderivative}
\[
f''(S^*) \ = \ 2n \cdot (S^*)^{-3} + \frac{1 - a}{a} \left( \frac{n + b}{n} \right)^{1/a} \cdot (S^*)^{(1 - 2a)/a}.
\]
To determine its sign, we reuse the condition \(f'(S^*) = 0\)
\[
f'(S^*) \ = \ -n \cdot (S^*)^{-2} + \left( \frac{n + b}{n} \right)^{1/a} \cdot (S^*)^{(1 - a)/a} \ = \ 0.
\]
Rewrite the equation by adding \(n \cdot (S^*)^{-2}\) to both sides and get
\[
\left( \frac{n + b}{n} \right)^{1/a} \cdot (S^*)^{(1 - a)/a} \ = \ n \cdot (S^*)^{-2}.
\]
Multiply both sides by \((S^*)^{-1}\), giving us
\[
\left( \frac{n + b}{n} \right)^{1/a} \cdot (S^*)^{(1 - 2a)/a} \ = \ n \cdot (S^*)^{-3}.
\]
Substituting the LHS of the previous equation into \(f''(S^*)\) simplifies the expression
\[
f''(S^*) \ = \ 2n \cdot (S^*)^{-3} + \frac{1 - a}{a} n \cdot (S^*)^{-3} \ = \ n \cdot \frac{1 + a}{a} \cdot (S^*)^{-3} \ > \ 0.
\]

From Equation \ref{f'(n)}, we evaluate the derivative at the endpoint \(S = n\)
\[
\begin{aligned}
f'(n) \ 
&= \ -n \cdot n^{-2} + \left( \frac{n + b}{n} \right)^{1/a} \cdot n^{(1 - a)/a} \ \\[6pt]
&= \ -n^{-1} + (n + b)^{1/a} n^{-1} \\[6pt]
&= \ \frac{(n + b)^{1/a} - 1}{n} \ > \ 0.
\end{aligned}
\]

We substitute the optimal value \(S^*\) into \(f(S)\) from Equation \ref{f(S*)}
\[
f(S^*) \ = \ \frac{n}{S^*} + a \left( \frac{n + b}{n} \right)^{1/a}  (S^*)^{1/a}.
\]
The first term is
\[
\begin{aligned}
\frac{n}{S^*} \ 
&= \ \frac{n}{n(n + b)^{-1/(a + 1)}} \ \\[6pt]
&= \ (n + b)^{1/(a + 1)}.
\end{aligned}
\]
The second term is
\[
\begin{aligned}
a \left( \frac{n + b}{n} \right)^{1/a}  (S^*)^{1/a} \ 
&= \ a \left( \frac{n + b}{n} \right)^{1/a} \left( n(n + b)^{-1/(a + 1)} \right)^{1/a} \\[6pt]
&= \ a (n + b)^{1/a} (n + b)^{-1/a(a + 1)} \\[6pt]
&= \ a (n + b)^{1/(a + 1)}.
\end{aligned}
\]
Add both terms together to obtain the final value
\[
\begin{aligned}
f(S^*) \ 
&= \ (n + b)^{1/(a + 1)} + a (n + b)^{1/(a + 1)} \\[6pt]
&= \ (a + 1) (n + b)^{1/(a + 1)}.
\end{aligned}
\]

\section{The Three-Dimensional Critical Point Problem}\label{appendix:3D-critical}

What about evolving the problem into a three-dimensional case? Similar to the two-dimensional problem, now consider a cuboid in three-dimensional space with coordinates from \( (0, 0, 0) \) to \( (L, M, N) \), inclusive. Without loss of generality, we assume \( L \geq M \geq N \). Then the target is to use the minimal number of egg drops to locate the critical point \( (l, m, n) \), where \(l\), \(m\), \(n\), \(L\), \(M\), and \(N\) are all integers, with \(l \in (0, L]\), \(m \in (0, M]\), and \(n \in (0, N]\), such that if a magical egg is dropped at a point \( (a, b, c) \), then
\begin{itemize}
\item if \( a < l \), \( b < m \), and \( c < n \), the egg remains intact, with no damage;
\item otherwise (i.e., if \( a \geq l \), \( b \geq m \), or \( c \geq n \)), the egg breaks.
\end{itemize}

\begin{theorem}\label{thm:3D}
If there are \(k\) eggs, then in the worst-case scenario the minimum number of drops in the worst case is
\[
P_3(k) \ \leq \ \lceil (k - 2) \cdot (L + M + N)^{1/(k - 2)} \rceil, \quad \text{for } k \geq 3,
\]
\end{theorem}

\subsection{Base Case: \( k = 3 \)}

It is impossible to determine all coordinates of the critical point using only one or two eggs, as at least three eggs are needed to guarantee the completion of searching along all three axes.

When three eggs are available, the only guaranteed strategy is to perform a linear search along each axis—using one egg per axis to determine each of the critical \( x \)-, \( y \)-, and \( z \)-coordinates. In the worst-case scenario, the eggs do not break until it is dropped at positions \( (L,0,0) \), \( (0,M,0) \) and \( (0,0,N) \), resulting in a total of \(L + M + N\) egg drops.

Note that when \( k = 3 \),
\[
P_3(3) \ \le \ \lceil (3 - 2) \cdot (L + M + N)^{1/(3 - 2)} \rceil = L + M + N,
\]
which exactly matches the number of drops required by the linear strategy, implying the base case holds at \(3\).

\subsection{Inductive Step}

Assume that with \( k \geq 3\) eggs, the minimum number of drops required in the worst case is
\[
P_3(k) \ \le \ \lceil (k - 2) \cdot (L + M + N)^{1/(k - 2)} \rceil.
\]
We aim to prove that, if \(P_3(k)\) is true for \(k\), then \(P_3(k + 1)\) is true, which indicates that, with \( k + 1 \) eggs, the worst-case minimum number of drops is
\[
P_3(k + 1) \ \le \ \lceil (k - 1) \cdot (L + M + N)^{1/(k - 1)} \rceil.
\]

\begin{figure}[h]
\begin{center}
	\scalebox{.22}{\includegraphics{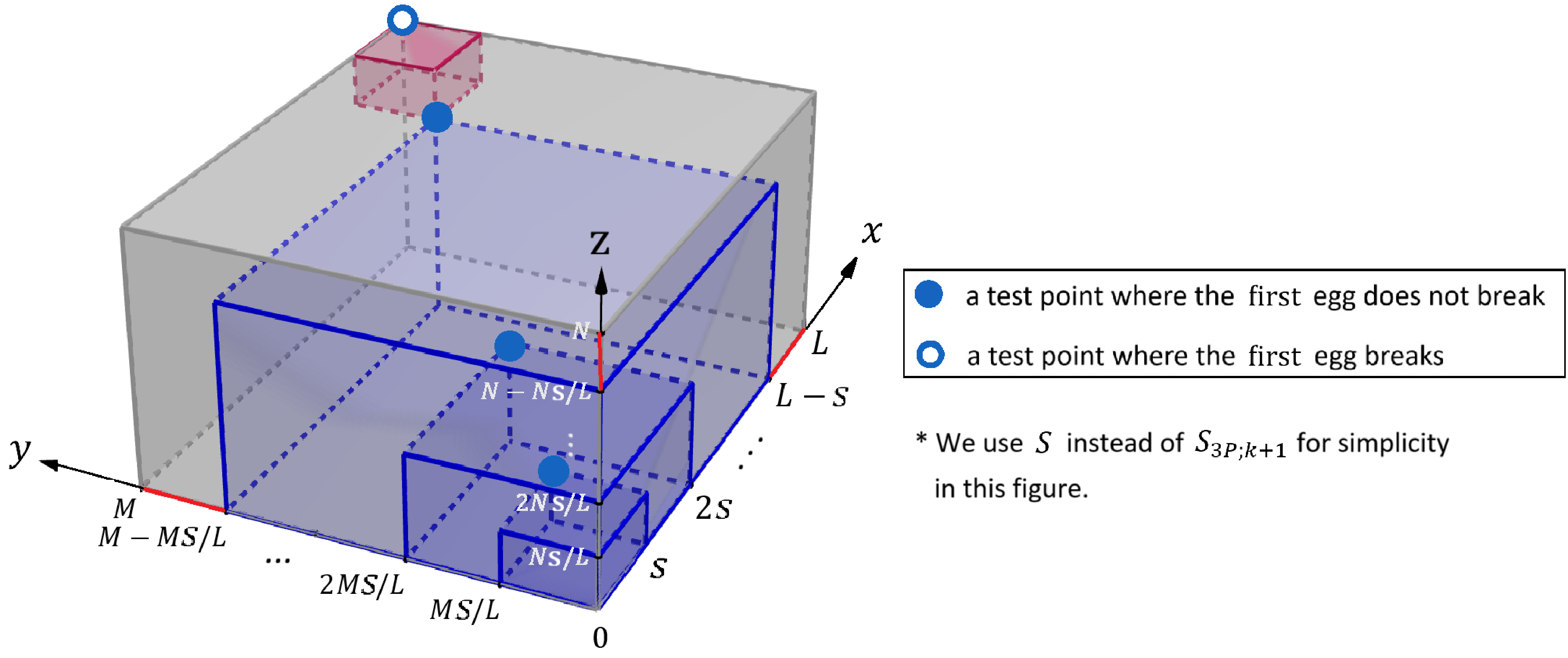}}
	\caption{\label{fig:Point3DWorstCase} The 3D critical point strategy.}
	\end{center}\end{figure}

We drop the first egg at points of the form
\((iS_{3P; k + 1}, iMS_{3P; k + 1} / L, iNS_{3P; k + 1}/L), \quad \text{for } i = 1, 2, 3, \ldots\), where \(S_{3P; k + 1}\) is a \(k\)-egg step size to be determined when searching for a \textit{critical point} (P) in the three-dimensional setting, until it breaks at \((L, M, N)\) (as shown in Figure \ref{fig:Point3DWorstCase}), which requires \( L/S_{3P; k + 1} \) drops in the worst case. Assuming the same convention as above, \(iMS_{3P; k + 1}/L\) and \(iNS_{3P; k + 1}/L\) can be any real numbers, not necessarily integers.

The critical point must lie within a sub-cuboid of dimensions \( S_{3P; k + 1} \), \( MS_{3P; k + 1}/L\), and \( NS_{3P; k + 1}/L \) shown in red. Since one egg has been used, we are now left with \( k \) eggs and must search this smaller red sub-cuboid under the same recursive strategy. By the inductive hypothesis, the number of drops required for this sub-problem is
\[
(k - 2) \cdot \left( S_{3P; k + 1} + \frac{MS_{3P; k + 1}}{L} + \frac{NS_{3P; k + 1}}{L} \right)^{1/(k - 2)} \ = \ (k - 2) \cdot \left( \frac{L + M + N}{L} \right)^{1/(k - 2)} \cdot S_{3P; k + 1}^{1/(k - 2)}.
\]

	Thus, the total number of drops in the worst case is diagonal tries plus linear tries, and can be expressed as the function
	\[
	f_{3P; k + 1}(S_{3P; k + 1}) \ \le \ \frac{L}{S_{3P; k + 1}} + (k - 2) \cdot \left( \frac{L + M + N}{L} \right)^{1/(k - 2)} \cdot S_{3P; k + 1}^{1/(k - 2)}.
	\]
	
	We now use Lemma \ref{lem:globalminab} as \(a = k - 2, b = M + N, \text{ and } n = L\) giving,
	\[
	S_{3P; k + 1}^* \ = \ L \cdot (L + M + N)^{-1/(k - 1)},
	\]
	and
	\[
	f_{3P; k + 1}(S_{3P; k + 1}^*) \ = \ (k - 1) \cdot (L + M + N)^{1/(k - 1)}.
	\]
	
	Thus, the inductive step holds: as long as given \(P_3(k)\) is true, with \( k + 1 \) eggs, the worst-case minimum number of drops is
	\[
	P_3(k + 1) \ \le \ \lceil (k - 1) \cdot (L + M + N)^{1/(k - 1)} \rceil,
	\]
	completing the proof. \(\hfill \Box\)

	\section{Analysis and Comparison} \label{analysis}
	
	Regarding the critical line \(x + y = V\) problem illustrated in Section \ref{critical line}, in order to learn which method is a better strategy, we perform only a quick comparison by keeping terms up to the second order in the Taylor expansions. We emphasize that the analysis in this section is not a rigorous proof. A fully rigorous error analysis remains a task for future work.
	
	We begin by define upper bound of the number of drops required under \textbf{Method One} as
	\[
	l_1(k) \ := \ k \cdot (M + N)^{1/k}, \{ k \in \mathbb{Z} \mid k \ge 2 \}.
	\]
	Similarly, define the upper bound of the number of drops required under \textbf{Method Two} as
	\[
	l_2(k) \ := \ (k - 1) \cdot M^{1/(k - 1)} + 1, \{ k \in \mathbb{Z} \mid k \ge 2 \}.
	\]
	We next examine the difference between these two bounds, which we believe are close to the true value. We define
	\begin{equation}\label{eq:l(k)}
l(k) \ := \ l_1(k) - l_2(k) \ = \ k \cdot (M + N)^{1/k} - \left( (k - 1) \cdot M^{1/(k - 1)} + 1\right) , \{ k \in \mathbb{Z} \mid k \ge 2 \}.
\end{equation}

The sign of \(l(k)\) determines which methods is more efficient.
\begin{itemize}

\item If \(l(k) < 0\), then \(l_1(k) < l_2(k)\), and \textbf{Method One} requires fewer drops.

\item If \(l(k) > 0\), then \(l_1(k) > l_2(k)\), and \textbf{Method Two} requires fewer drops.

\item The magnitude \( \left|l(k) \right| \) reflects how much better one method performs compared with the other.
\end{itemize}

\begin{figure}[h]
\centering
\scalebox{.2}{\includegraphics{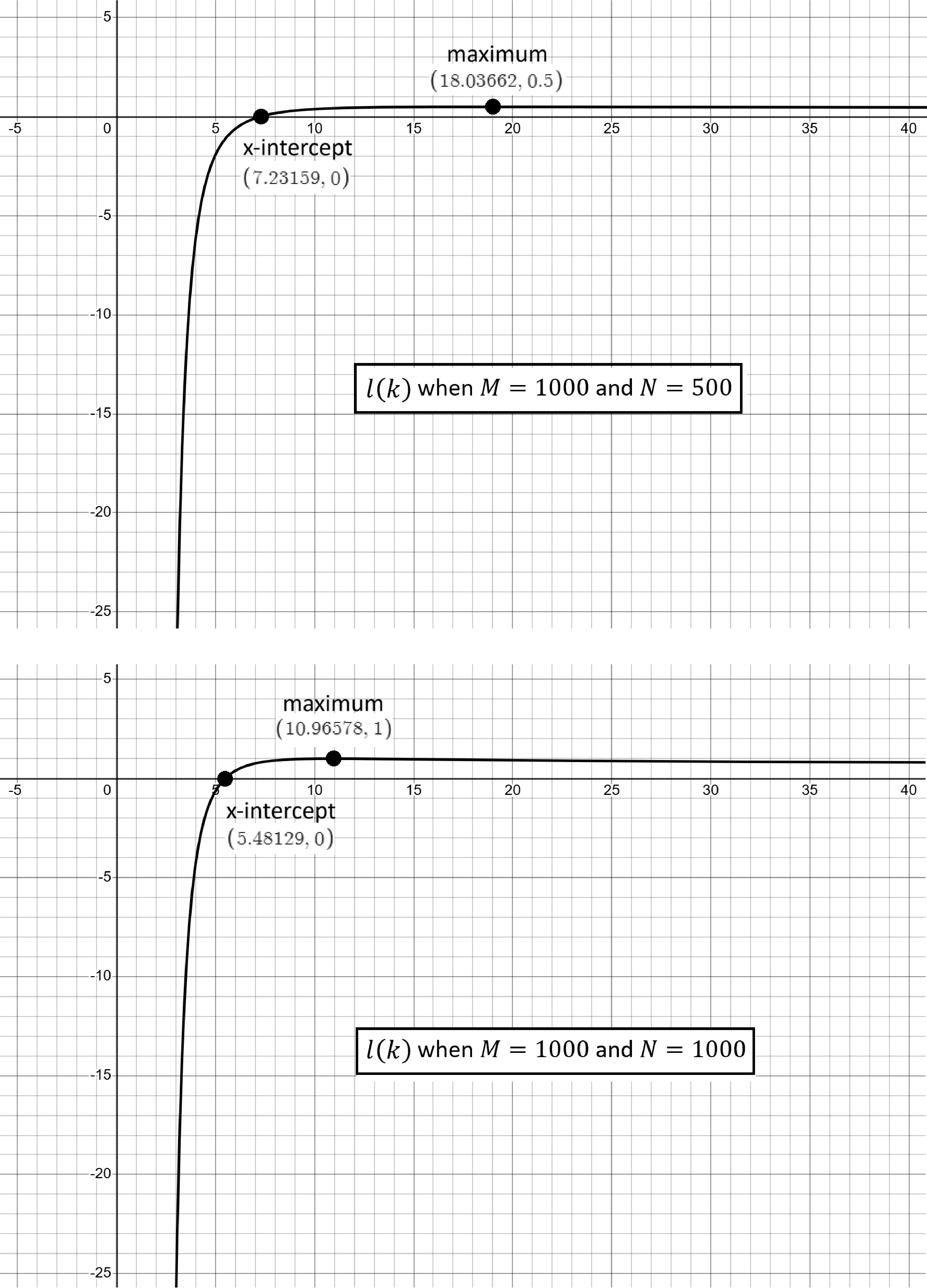}}
\caption{\label{fig:LineXYanalysis} Sample plots showing the behavior of \(l(k)\).}
\end{figure}

From numeral observations as shown in Figure \ref{fig:LineXYanalysis}, we find that for small values of \(k\) (for example, when only \(2\) or \(3\) eggs are available), \textbf{Method One} is significantly better than Method Two. For larger values of \(k\) (i.e., when we have a larger supply of eggs, such as \(10\) or \(20\)), Method Two performs only slightly (and negligibly) better than Method One. So, overall, \textbf{Method One} provides  a better strategy.

\begin{proof}[Approximate]
To understand why \(l(k)\) behaves as shown in Figure \ref{fig:LineXYanalysis}, we use the Taylor series expansion\footnote{As several of us are in pre-highschool, we used ChatGPT for a suggestion and from that we learnt that we could use the Taylor expansion.} at arbitrarily large \(k\) to get
\[
(M + N)^{1/k} \ = \ e^{1/k \cdot \ln(M + N)} \ = \ \sum_{n = 0}^{\infty} \frac{(\frac{1}{k} \ln(M + N))^n}{n!} \ = \ \sum_{n = 0}^{\infty} \frac{(\ln(M + N))^n}{k^n \cdot n!}
\]
and
\[
M^{1/(k - 1)} \ = \ e^{1/(k -1) \cdot  \ln M} \ = \ \sum_{n = 0}^{\infty} \frac{(\frac{1}{k - 1} \ln M)^n}{n!} \ = \ \sum_{n = 0}^{\infty} \frac{(\ln M)^n}{(k - 1)^n \cdot n!}.
\]
Substituting into Function \ref{eq:l(k)} gives
\[
l(k) \ = \ k \cdot \sum_{n = 0}^{\infty} \frac{(\ln(M + N))^n}{k^n \cdot n!} - \left( (k - 1) \cdot \sum_{n = 0}^{\infty} \frac{(\ln M)^n}{(k - 1)^n \cdot n!} \right) - 1.
\]

Let \(n\) be the number of correction terms we keep and \(T_n\) be the \(n\)th-degree Taylor polynomial of \(l(k)\).
When using a Taylor polynomial \(T_n\) to approximate a function, we need to decide how large to take \(n\) to ensure that we achieve a desired accuracy \cite{Stewart}.
In order to do so, we graph \(l(k)\) and \(T_n(k)\) where \(n = 1, 2, 3, \text{ and } 4\) as shown in Figure \ref{fig:LineXYerror}.

Our goal is to determine the sign of \(l(k)\) --- whether it is positive or negative --- and \(T_2(k)\) can well reflect the feature of \(l(k)\), so keeping the first two correction terms is sufficient to achieve this. We only carry out a second-order analysis here; one could make this rigorous by a more careful treatment of the error terms.

\begin{figure}[h]
	\centering
	\scalebox{.2}{\includegraphics{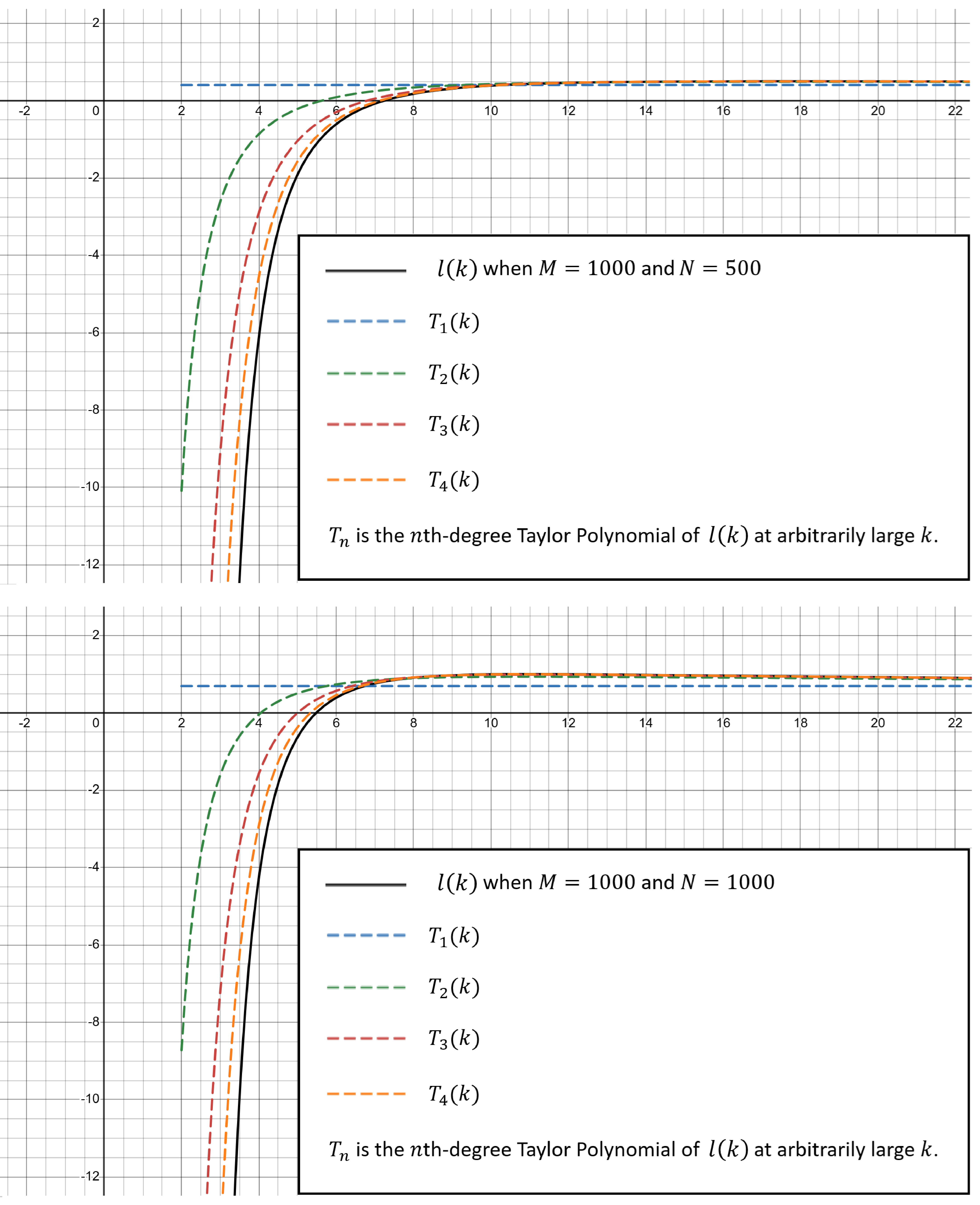}}
	\caption{\label{fig:LineXYerror} Error estimation of \(l(k)\).}
\end{figure}

Keeping two correction terms in \(l(k)\), we get
\[
T_2(k) \ = \ k\left( 1 + \frac{\ln(M + N)}{k} + \frac{(\ln(M + N))^2}{2k^2}\right)  - (k - 1)\left( 1 + \frac{\ln M}{k - 1} + \frac{(\ln M)^2}{2(k - 1)^2} \right) - 1.
\]
Simplifying gives

\[
T_2(k) \ = \ \ln \left( 1 + \frac{N}{M}\right)  + \frac{(\ln(M + N))^2}{2k} - \frac{(\ln M)^2}{2(k - 1)}.
\]

\textbf{The leading term}, \( \ln (1 + N/M)\), lies between \(0\) and \(\ln 2 \approx 0.693\).
In \textbf{the correction terms}, the difference between \((\ln(M + N))^2 \text{ and } (\ln M)^2\) is small, since the difference between \(\ln(M + N) \text{ and } \ln M\) is \( \ln (1 + N/M)\), which is smaller than \(\ln 2 \approx 0.693\). The difference between \(2k\) and \(2(k - 1)\) is \(2\).

For \textbf{small values of \(k\)}, the relative difference between \(2k\) and \(2(k - 1)\) is large, implying that \((\ln(M + N))^2/2k\) is much smaller than \((\ln M)^2/2(k - 1)\). Thus, the negative correction term \((\ln M)^2/2(k - 1)\) domains, making \(l(k)\) significantly negative.
For \textbf{large values of \(k\)}, \(2k\) and \(2(k - 1)\) are nearly equal, implying that the correction terms almost cancel out, and leaving \(l(k) \approx \ln (1 + N/M)\), a small positive number in the range of \((0,1)\).

As \(N\) increases (or approaches to \(M\)), the positive terms of \(T_2(k)\) grow larger, causing the function to become positive more rapidly and shifting its \(x\)-intercept to smaller values of \(k\). However, this change does not affect the overall qualitative features of the graph.

In summary, for small \(k\) (for example, if \(2\) or \(3\) eggs are available), \(\textbf{Method One}\) outperforms Method Two. For large \(k\) (such as \(k = 10\) or \(k = 20\)), Method Two is slightly better, but the advantage is negligible. Overall, \(\textbf{Method One}\) is the more efficient strategy.
\end{proof}



\begin{thebibliography}{10990}

\bibitem[Bo]{Boardman}
\newblock Michael Boardman.
\textit{Egg Drop Numbers}, Mathematics Magazine \textbf{77} (2004), no. 5. 368-372.

\bibitem[Br]{Brilliant}
\newblock Brilliant. \bburl{https://brilliant.org/wiki/egg-dropping/}.

\bibitem[Fr]{Fr}
\newblock Joel Franklin.
\emph{Mathematical Methods of Economics: Linear and Nonlinear Programming, Fixed-Point Theorem}, Springer-Verlag, New York, 1980.

\bibitem[GG]{Geeks}
\newblock Geeks for Geeks. \bburl{https://www.geeksforgeeks.org/dsa/egg-dropping-puzzle-dp-11/}.

\bibitem[KVW]{KonhauserVellemanWagon}
\newblock J. D. E. Konhauser, D. Velleman, and S. Wagon.
\textit{Which Way Did the Bicycle Go?: And Other Intriguing Mathematical Mysteries}, Dolciani Mathematical Expositions. Mathematical Association of America. \textbf{18} (1996), no. 18.

\bibitem[Mi1]{Miller1}
\newblock Steven J. Miller.
\textit{Mathematics of Optimization: How to do things faster}, AMS, Pure and Applied Undergraduate Texts \textbf{30} (2017), Chapter 9.

\bibitem[Mi2]{Miller2}
\newblock Steven J. Miller, Daniel Turek, and Haoyu Sheng.
\emph{When Rooks Miss: Probability through Chess}, the College Mathematics Journal \textbf{52} (2021), no. 2, 82--93.


\bibitem[Mi3]{Miller3}
\newblock Steven J. Miller.
\emph{Math 331: The Little Questions (Fall 2024), Lecture 9}.
Course homepage:	\url{https://web.williams.edu/Mathematics/sjmiller/public_html/331Fa24/}.
Video: \url{https://youtu.be/pGixuGkTuQM}.
Slides (starting on page 124):	\url{https://web.williams.edu/Mathematics/sjmiller/public_html/331Fa24/math331fa24slides.pdf}.


\bibitem[Mo]{SM}
\newblock Spencer Mortensen. \bburl{https://spencermortensen.com/articles/egg-problem/}.

\bibitem[Na]{Nair}
\newblock Umesh Nair.
\textit{Analysis of Puzzles}, \textbf{2} (2006).

\bibitem[PW]{ParksWills}
\newblock Harold R. Parks and Dean C. Wills.
\textit{Two Eggs Any Style Generalizing Egg Drop Experiments}, Recreational Mathematics Magazine. \textbf{12} (2025), no. 20. 1-18.

\bibitem[Sn]{Sniedovich}
\newblock Moshe Sniedovich.
\textit{OR/MS Games: 4. The Joy of Egg-Dropping in Braunschweig and Hong Kong}, \textbf{4} (2003), no. 1. 48-64.

\bibitem[St]{Stewart}
\newblock James Stewart.
\textit{Calculus: Early Transcendentals}, Cengage Learning. Chapter 11.

\bibitem[Yo]{Youtube}
\newblock Simply Logical, Youtube
\textit{Egg Dropping Puzzle with 2 Eggs and 100 Floors -- Microsoft Interview Puzzles}
\bburl{https://www.youtube.com/watch?v=uBhSIKLlvdk}.
\end{thebibliography}
\end{document}